\documentclass[aap,MSNbibl,citesort,dvips]{arximspdf}
\usepackage{mathbh}
\usepackage{graphicx}

\doi{10.1214/10-AAP728}
\volume{21}
\issue{4}
\pubyear{2011}
\firstpage{1400}
\lastpage{1435}

\makeatletter
\newtheorem{claim}[tthm]{Claim}

\newproclaim{re}{Remark}

\newcommand{\eqref}[1]{(\ref{#1})}
\newtheorem{tthm}{Theorem}[section]
\newtheorem{lmm}[tthm]{Lemma}

\newtheorem{prop}[tthm]{Proposition}

\makeatother

\begin{document}
\begin{frontmatter}

\title{Random graphs with a given degree sequence}
\runtitle{Random graphs with a given degree sequence}

\begin{aug}
\author[A]{\fnms{Sourav} \snm{Chatterjee}\corref{}\ead[label=e1]{sourav@cims.nyu.edu}\thanksref{t1}},
\author[B]{\fnms{Persi} \snm{Diaconis}\ead[label=e2]{diaconis@math.stanford.edu}}
\and
\author[C]{\fnms{Allan} \snm{Sly}\ead[label=e3]{allansly@microsoft.com}}

\thankstext{t1}{Supported in part by NSF Grant DMS-07-07054 and a Sloan
Research Fellowship.}
\runauthor{S. Chatterjee, P. Diaconis and A. Sly}
\affiliation{New York University, Stanford University and Microsoft Research}
\address[A]{S. Chatterjee\\
Courant Institute of Mathematical Sciences\\
New York University\\
251 Mercer Street\\
New York, New York 10012\\
USA\\
\printead{e1}} %adresu isvedimo komanda gale!
\address[B]{P. Diaconis\\
Department of Mathematics\\
Stanford University\\
Stanford, California 94305\\
USA\\
\printead{e2}}
\address[C]{A. Sly\\
Theory Group\\
Microsoft Research, One Microsoft Way\\
Redmond, Washington 98052\\
USA\\
\printead{e3}}
\end{aug}

% HISTORY:
\received{\smonth{5} \syear{2010}}
\revised{\smonth{8} \syear{2010}}

% ABSTRACT
%
\begin{abstract}
Large graphs are sometimes studied through their degree sequen\-ces
(power law or regular graphs). We study graphs that are uniformly
chosen with a given degree sequence. Under mild conditions, it is shown
that sequences of such graphs have graph limits in the sense of
Lov\'asz and Szegedy with identifiable limits. This allows simple
determination of other features such as the number of triangles. The
argument proceeds by studying a natural exponential model having the
degree sequence as a sufficient statistic. The maximum likelihood
estimate (MLE) of the parameters is shown to be unique and consistent
with high probability. Thus $n$ parameters can be consistently
estimated based on a sample of size one. A fast, provably convergent,
algorithm for the MLE is derived. These ingredients combine to prove
the graph limit theorem. Along the way, a continuous version of the
Erd\H{o}s--Gallai characterization of degree sequences is derived.
\end{abstract}

% KEYWORDS
%

\end{frontmatter}

%s1 ###
\section{Introduction}\label{intro}
%s1.1 ###
\subsection{Graphs with a given degree sequence}

Let $G$ be an undirected simple graph on $n$ vertices and let
$d_1,\ldots, d_n$ be the degrees of the vertices of~$G$. The vector
$\mathbf{d}:= (d_1,\ldots, d_n)$ is usually called the \textit{degree
sequence} of~$G$. Correspondingly, the \textit{degree distribution} of
$G$ is the probability distribution function $F$ supported on $[0,1]$,
defined as
\[
F(x) := \frac{|\{i\dvtx d_i \le nx\}|}{n}.
\]
In other words, if a vertex is chosen uniformly at random, then the
degree of that vertex, divided by $n$, is a random variable with
probability distribution function $F$.

In recent years, the degree distributions of real world networks have
received wide attention. The surveys \cite{newman03,newmanetal06}
contain many references as does the detailed account in
\cite{blitzsteindiaconis09}. The enthusiasm of some authors for ``scale
free'' or ``power law graphs'' has also generated much controversy
\cite{lietal05,willingeretal09} which serves as additional motivation
for the present paper.
% (although we are concerned with dense graphs, whereas power law
%graphs are typically sparse; however, with some effort, it may be
%possible to extend our approach to the sparse case in a future
%project).

The interest in degree distributions stems from the fact that the
degree sequences of real world networks sometimes appear to have power
law behavior that is very different than those occurring in classical
models of random graphs, like the Erd\H{o}s--R\'enyi model
\cite{erdosrenyi60}. Researchers have tried various ways of
circumventing this problem. An obvious solution is to build random
graph models that are forced to give us the degree distribution that we
want and then deduce other features by simulation or mathematics. A
natural way to do this is to choose a graph uniformly at random from
the set of all graphs with a given degree sequence. One frequent
appearance of this model is for random regular graphs \cite{wormald99}.
As explained in \cite{blitzsteindiaconis09}, Section 13, the model also
arises in testing if the exponential family with degree sequence as
sufficient statistic fits a given data set. See \cite{tsourakakis08}
for applications where the number of triangles is wanted. The paper
\cite{blitzsteindiaconis09} has useful ways of simulating graphs with a
given degree sequence and an extensive survey of the (mostly
nonrigorous) literature for this model. Some rigorous results are also
available in the ``sparse case,'' for example, those
in~\cite{molloyreed95,molloyreed98}.

At this point, a gap between our motivation and our theory must be
pointed out: the present paper deals with \textit{dense} graphs with a
given degree sequence (roughly, graphs whose number of edges is
comparable to the square of the number of vertices), whereas much of
the literature cited above, for example, power law graphs, revolves
around sparse graphs. As of now, our theorems are not directly
applicable in the sparse setting, although there is certainly hope for
future progress.

%In recent work Barvinok and Hartigan~\cite{barvinokHartigan10} under
%mild regularity conditions have given asymptotic formulas for the
%number of graphs with a given degree sequence. (See
%well-known counting problem.)

In a recent series of papers
\cite
{barvinok08,barvinok10,barvinokHartigan09b,barvinokHartigan09c,barvinokHartigan09a,barvinokHartigan10},
Barvinok and Hartigan have looked at problems related to the structure
of directed and undirected (dense) graphs with given degree sequence.
The Barvinok and Hartigan work, especially~\cite{barvinokHartigan10},
is related to the present paper. This is explained at the end of this
Introduction after we have stated our main theorems.

One of the objectives of this article is to give a rather precise
description of the structure of random (dense) graphs with a given
degree sequence via the notion of graph limits introduced recently by
Lov\'asz and Szegedy~\cite{lovaszszegedy06} and developed by Borgs
et~al.~\cite{borgsetal06,borgsetal08,borgsetal07}. See also the related work of
Diaconis and Janson \cite{diaconisjanson08} and Austin \cite{austin08}
which traces this back to work of Aldous \cite{aldous81} and
Hoover~\cite{hoover82}. This gives, in particular, a way to write down
exact formulas for the expected number of subgraphs of a given type
without simulation.

Before stating our result, we need to introduce the notion of graph
limits. We quote the definition verbatim from \cite{lovaszszegedy06}
(see also \cite{borgsetal08,borgsetal07,diaconisjanson08}). Let~$G_n$
be a sequence of simple graphs whose number of nodes tends to infinity.
For every fixed simple graph $H$, let $|\hom(H, G)|$ denote the number
of homomorphisms of $H$ into $G$ [i.e., edge-preserving maps $V(H)
\rightarrow
V(G)$, where~$V(H)$ and $V(G)$ are the vertex sets]. This number is
normalized to get the homomorphism density
%
%e1 ###
\begin{equation}\label{homdens}
t(H,G) := \frac{|\hom(H, G)|}{|V(G)|^{|V(H)|}}.
\end{equation}
This gives the probability that a random mapping $V(H) \rightarrow
V(G)$ is a
homomorphism.

Suppose that the graphs $G_n$ become more and more similar in the sense
that $t(H, G_n)$ tends to a limit $t(H)$ for every $H$. One way to
define a limit of the sequence $\{G_n\}$ is to define an appropriate
limit object from which the values $t(H)$ can be read off.

The main result of \cite{lovaszszegedy06} (following the earlier
equivalent work of Aldous \cite{aldous81} and Hoover \cite{hoover82})
is that indeed there is a natural ``limit object'' in the form of a
symmetric measurable function $W \dvtx [0, 1]^2 \rightarrow[0, 1]$ [we
call $W$ symmetric if $W(x, y) = W (y, x)$]. Conversely, every such
function arises as the limit of an appropriate graph sequence. This
limit object determines all the limits of subgraph densities: if $H$ is
a simple graph with $V(H) = [k] = \{1, \ldots, k\}$, then
\[
t(H,W) = \int_{[0,1]^k}\prod_{(i,j)\in E(H)} W(x_i, x_j)\,dx_1\,\cdots\,dx_k.
\]
Here $E(H)$ denotes the edge set of $H$.

Intuitively, the interval $[0,1]$ represents a ``continuum'' of
vertices and~$W(x,\break y)$ denotes the probability of putting an edge
between $x$ and $y$. For example, for the Erd\H{o}s--R\'enyi graph
$G_{n,p}$, if $p$ is fixed and $n \rightarrow\infty$, then~the limit
graph is
represented by the function that is identically equal to $p$
on~$[0,1]^2$.

Convergence of a sequence of graphs to a limit has many consequences.
From the definition, the count of fixed size subgraphs converges to the
right- hand side of the expression for $t(H,W)$ given above. More
global parameters also converge. For example, the degree distribution
converges to the law of $\int_0^1 W(U,y)\,dy$ where $U$ is a random
variable distributed uniformly on $[0,1]$. Similarly, the distribution
function of the eigenvalues of the adjacency matrix converges. More
generally, a graph parameter is a function from the space~of graphs
into a space $\mathcal{X}$ which is invariant under isomorphisms. If
$\mathcal{X}$~is a~to\-pological space, we may ask which graph parameters
are continuous with respect to the topology induced by graph limits.
This is called ``property~tes\-ting'' in the computer science theory
literature which has identified many continuous graph parameters. See
the surveys \cite{austintao10,borgsetal08} for pointers to the~lite\-rature.

We are now ready to state our result about the limit of graphs with~%
given degree sequences. Suppose that for each $n$, a degree sequence
$\mathbf{d}^n = (d_1^n, \ldots, d^n_n)$ is given. Without loss of generality,
assume that $d_1^n\ge d_2^n \ge\cdots\ge d_n^n$. We say\vadjust{\eject} that the
sequence $\{\mathbf{d}^n\}$ has a \textit{scaling limit} if there is a
nonincreasing function~$f$ on $[0,1]$ such that
%
%e2 ###
\begin{equation}\label{degconv}
\lim_{n\rightarrow\infty}\Biggl(\biggl| \frac{d_1^n}{n} - f(0) \biggr| + \biggl| \frac
{d_n^n}{n} - f(1)
\biggr| + \frac1n\sum_{i=1}^n \biggl|\frac{d_i^n}{n} - f\biggl(\frac{i}{n}\biggr)\biggr| \Biggr) = 0.
\end{equation}
It is not difficult to prove by a simple compactness argument that any
sequence $\{\mathbf{d}^n\}$ of degree sequences has a subsequence that
converges to a~scaling limit in the above sense. %Therefore, one cannot
%have a `more general' notion of scaling limits of degree sequences.
Note that convergence in the above sense can be stated equivalently in
terms of convergence of degree distributions: $d_1^n/n \rightarrow f(0)$,
$d_n^n/n \rightarrow f(1)$ and $D_n/n \rightarrow f(U)$ in
distribution, where $D_n$ is
a~randomly (uniformly) chosen $d_i^n$ and $U$ is uniformly distributed
on $[0,1]$.

The need to control $d_1^n$ and $d_n^n$ arises from the need to
eliminate ``outlier'' vertices that connect to too many or too few
nodes which takes the degree sequence too close to the
Erd\H{o}s--Gallai boundary (see below). Since $d_i^n$ is decreasing in
$i$, outliers can be eliminated by simply controlling $d_1^n$ and
$d_n^n$. The need to eliminate outliers, on the other hand, arises from
technical aspects of our analysis.

Define ${D'[0,1]}$ to be the set of nonincreasing functions on $[0,1]$
which are left continuous on $(0,1)$. The reason for imposing
left-continuity is the following: When the scaling limit of a degree
sequence is discontinuous, it is not uniquely defined but there always
exists a unique limit in ${D'[0,1]}$. We could have as well chosen
right-continuous. %So we shall restrict our attention to limits in $

For each $n$, let $G_n$ be a random graph chosen uniformly from the set
of all simple graphs with degree sequence $\mathbf{d}^n$. Let $f$ be the
scaling limit of the sequence $\{\mathbf{d}^n\}$ in the sense defined above.
Our objective is to compute the limit of the sequence $\{G_n\}$ in
terms of the scaling limit of $\mathbf{d}^n$. We endow the set of scaling
limits (i.e., ${D'[0,1]}$) with the topology induced by a modified~$L^1$
norm $\| \cdot\|_{{1'}}$ given by
\[
\| f \|_{{1'}} := |f(0)| +|f(1)| + \int_0^1 |f(x)|\,dx.
\]
The choice of this norm is necessitated by the need to make it
compatible with our previous notion of convergence of degree sequences.
%%the need to control the behavior of the largest and smallest degrees.
%In particular it will avoid the case that $d_1^n=n$ or $d_n^n=0$ where
%the MLE described in the next section is not defined.

Not all functions can be scaling limits of degree sequences. Let
$\mathcal{F}$
be the set of functions in ${D'[0,1]}$ that can be obtained as scaling
limits of degree sequences in the sense stated above. By a simple
diagonal argument, it is easy to see that $\mathcal{F}$ is a closed
subset of
${D'[0,1]}$ under the topology of the modified $L^1$ norm. It is shown
in Proposition \ref{egcond} that $\mathcal{F}$ has nonempty interior.

\begin{tthm}\label{degthm}
Let $G_n$ and $f$ be as above. Suppose that $f$ belongs to the
topological interior of the set $\mathcal{F}$ defined above. Then
there exists
a unique function $g\dvtx [0,1]\rightarrow\mathbb{R}$ in ${D'[0,1]}$ such
that the function
\[
W(x,y) := \frac{e^{g(x)+g(y)}}{1+e^{g(x)+g(y)}}
\]
satisfies, for all $x\in[0,1]$,
\[
f(x) = \int_0^1 W(x,y)\,dy.
\]
In this situation, the sequence $\{G_n\}$ converges almost surely to
the limit graph represented by the function $W$.
\end{tthm}

Theorem 1.1 can be useful only if we can provide a simple way of
checking whether $f$ belongs to the interior of $\mathcal{F}$. (Being
the limit
of a sequence of degree sequences, it is clear that $f\in\mathcal
{F}$. The
nontrivial question is whether~$f$ is in the interior.) The following
result gives an easily verifiable equivalent condition.

\begin{prop}\label{egcond}
A function $f\dvtx [0,1]\rightarrow[0,1]$ in ${D'[0,1]}$ belongs to the
interior of
$\mathcal{F}$ if and only if:
\begin{enumerate}[(ii)]
\item[(i)] there are two constants $c_1 > 0$ and $c_2< 1$ such that
$c_1\le f(x)\le c_2$ for all $x\in[0,1]$ and
\item[(ii)] for each
$x\in(0,1]$,\vspace*{1pt}
\[%\label{eg}
\int_x^1 \min\{f(y), x\}\,dy + x^2 - \int_0^x f(y)\,dy > 0.\vspace*{1pt}
\]
\end{enumerate}
\end{prop}

\begin{re}
Condition (ii) in the above result is a continuum
version of the well-known Erd\H{o}s--Gallai criterion
\cite{erdosgallai60}: Suppose $d_1 \ge d_2 \ge\cdots\ge d_n$ are
nonnegative integers. The Erd\H{o}s--Gallai criterion says that
$d_1,\ldots,d_n$ can be the degree sequence of a simple graph on $n$
vertices if and only if $\sum_{i=1}^n d_i$ is even and for each $1\le
k\le n$,\vspace*{1pt}
\[
\sum_{i=1}^k d_i \le k(k-1) + \sum_{i=k+1}^n \min\{d_i, k\}.\vspace*{1pt}
\]
(See \cite{mahadevpeled95} for extensive discussions
and eight equivalent conditions.) %The choice of the modified $L^1$
%norm is made to ensure that
\end{re}

\begin{re}
When the scaling limit $f$ is continuous,
convergence in the modified $L^1$ norm is the same as supnorm
convergence. In particular, for continuous scaling limits
Theorem~\ref{degthm} and Proposition~\ref{egcond} both hold if we
replace ${D'[0,1]}$ with $C[0,1]$ and redefine $\mathcal{F}$
analogously under
the supnorm topology.
\end{re}

\begin{re}
As an example, consider the limit of the
Erd\H{o}s--R\'enyi graph $G(n,p)$ as $n\rightarrow\infty$. Here
$f(x) = p$ for
all $x$. Condition (ii) becomes $(1-x) \min\{p,x\} + x^2 - px > 0$
for all $x$. Considering the two cases $x\ge p$ and $x< p$ it is easy
to see that this holds, so Erd\H{o}s--R\'enyi graphs are in the
interior of $\mathcal{F}$ for any fixed $p$, $0 < p < 1$.
\end{re}

\begin{re}
In a recent article \cite{mckay10} (following up on
the older work \cite{mckaywormald90}), McKay has computed subgraph
counts in random graphs with a given degree sequence. However, McKay's
results hold only if either the graph is sparse or the graph is dense
but all degrees are within $n^{1/2+\varepsilon}$ of the average degree. Thus,
it may be possible to recover Theorem \ref{degthm} from McKay's results
when the limit shape is a constant function but not in other cases.

The next natural question is whether one can feasibly compute the
function $g$ in Theorem \ref{degthm} for a given $f$. It turns out that
this is a central issue in the whole analysis. In fact, to prove
Theorem \ref{degthm} we analyze a related statistical model;
computation of the maximum likelihood estimate in that model leads to
an algorithm for computing $g$ which, in turn, yields a proof of
Theorem \ref{degthm}. The statistical model is discussed next.
\end{re}

%s1.2 ###
\subsection{Statistics with degree sequences}\label{statsec}
Informally, if the degree sequence captures the information in a graph,
different graphs with the same degree sequence are judged equally
likely. This can be formalized by saying that the degree sequence is a
sufficient statistic for a probability distribution on graphs. The
Koopman--Pitman--Darmois theorem forces this distribution to be of
exponential form. This approach to model building is explained and
developed in \cite{lauritzen88}. Diaconis and Freedman
\cite{diaconisfreedman84} give a version of the
Koopman--Pitman--Darmois theorem for discrete exponential families. The
approach is also standard fare in statistical mechanics where the
uniform distribution on graphs with fixed degree sequence is called
``micro-canonical'' and the exponential distribution is called
``canonical'' (see \cite{parknewman04}). It turns out
that the exponential model has a simple description in terms of
independent Bernoulli random variables.

Given a vector $\bolds{\beta}= (\beta_1,\ldots,\beta_n)\in\mathbb
{R}^n$, let
$\mathbb{P}_{\bolds{\beta}}$ be the law of the undirected random graph
on $n$ vertices
defined as follows: for each $1\le i\ne j\le n$, put an edge between
the vertices $i$ and $j$ with probability
\[
p_{ij} := \frac{e^{\beta_i + \beta_j}}{1+e^{\beta_i + \beta_j}}
\]
independently of all other edges. Thus, if $G$ is a graph with degree
sequence $d_1,\ldots,d_n$, the probability of observing $G$ under
$\mathbb{P}_{\bolds{\beta}}$ is
\[
\frac{e^{\sum_i \beta_i d_i}}{\prod_{i<j} (1+e^{\beta_i +\beta_j})}.
\]
Henceforth, this model of random graphs is called the ``$\bolds{\beta
}$-model.''
This model was considered by Holland and Lienhardt
\cite{hollandleinhardt81} in the directed case and by Park and Newman
\cite{parknewman04} and Blitzstein and Diaconis
\cite{blitzsteindiaconis09} in the undirected case. It is a close
cousin to the Bradley--Terry model for rankings [which itself goes back
(at least) to Zermelo]. See \cite{hunter04} for extensive
references. The $\bolds{\beta}$-model is also a simple version of a
host of
exponential models actively in use for analyzing network data. We will
not try to survey this vast literature but recommend the extensive
treatments in \cite{newman03,jackson08,robinsetal07}. The website for the International
Network for Social Network Analysis contains further information.

Suppose a random graph $G$ is generated from the $\bolds{\beta
}$-model where
$\bolds{\beta}\in\mathbb{R}^n$ is unknown. Is it possible to
estimate $\bolds{\beta}$ from the
observed $G$? It is not difficult to show that the maximum likelihood
estimate (MLE) $\hat{\bolds{\beta}}$ of $\bolds{\beta}$ must
satisfy the system of
equations
%
%e3 ###
\begin{equation}\label{mle}
d_i = \sum_{j\ne i} \frac{e^{\hat{\beta}_i+\hat{\beta
}_j}}{1+e^{\hat{\beta}_i +\hat{\beta}_j}},\qquad
i=1,\ldots,n,
\end{equation}
where $d_1,\ldots,d_n$ are the degrees of the vertices in the observed
graph $G$. Questions may arise about the existence, uniqueness and
accuracy of the MLE. Since the dimension of the parameter space grows
with $n$, it is not clear if this is a ``good'' estimate of $\bolds
{\beta}$ in
the traditional sense of consistency in statistical estimation theory.

The following theorem shows that under certain mild assumptions on
$\bolds{\beta}$, there is a high chance that the MLE exists, is
unique and
estimates $\bolds{\beta}$ with uniform accuracy in all coordinates.

\begin{tthm}\label{mainthm}
Let $G$ be drawn from the probability measure $\mathbb{P}_{\bolds
{\beta}}$ and let
$d_1,\ldots, d_n$ be the degree sequence of $G$. Let $L := \max_{1\le
i\le n} |\beta_i|$. Then there is a constant $C(L)$ depending only on
$L$ such that with probability at least $1- C(L)n^{-2}$, there exists a
unique solution $\hat{\bolds{\beta}}$ of the maximum likelihood equations
\eqref{mle}, that satisfies
\[
\max_{1\le i\le n}|\hat{\beta}_i -\beta_i|\le C(L)\sqrt{\frac
{\log n}{n}}.
\]
\end{tthm}

It may seem surprising that all $n$ parameters can be accurately
estimated from a single realization of the graph. However, one needs to
observe that there are, in fact, $n(n-1)/2$ independent random
variables lurking in the background (namely, the indicators whether
edges are present or not). There is a well-known heuristic that in a
$p$-parameter model with $m$ observations, ``the usual asymptotics''
work provided that $p^2/m$ tends to zero as $m$ tends to infinity. See
\cite{portnoy84,portnoy85,portnoy88,portnoy91} for details (and
counter examples). In our model $p=n$ and $m= n(n-1)/2$, so $p^2/m$
does not tend to zero but stays bounded. The heuristic, although not
directly applicable, hints at a reason why one can expect estimability
of parameters.

In work closer to the present paper, Simons and Yao \cite{simonsyao99}
studied the Bradley--Terry model for comparing $n$ contestants. Here a
random orientation of the complete graph on $n$ vertices is chosen
based on ``player $a$ beats player $b$ with probability
$\theta(a)/[\theta(a) + \theta(b)]$.'' They show that MLE is consistent
here as well. Hunter \cite{hunter04} shows that the MM algorithm also
behaves well in this problem.

The next theorem characterizes all possible expected degree sequences
of the $\bolds{\beta}$-model as $\bolds{\beta}$ ranges over
$\mathbb{R}^n$. The nice feature is
that no degree sequence is left out.

\begin{tthm}\label{meanspace}
Let $\mathcal{R}$ denote the set of all expected degree sequences of random
graphs following the law $\mathbb{P}_{\bolds{\beta}}$ as $\bolds
{\beta}$ ranges over $\mathbb{R}^n$. Let
$\mathcal{D}$ denote the set of all possible degree sequences of undirected
graphs on $n$ vertices. Then
\[
\operatorname{conv}(\mathcal{D}) = \overline{\mathcal{R}},
\]
where $\operatorname{conv}(\mathcal{D})$ denotes the convex hull of
$\mathcal{D}$ and $\overline
{\mathcal{R}}$ is the topological closure of~$\mathcal{R}$. %Taking
%closure is
%necessary because $\mr$ is not a closed set.
\end{tthm}

Incidentally, the convex hull of $\mathcal{D}$ is a well-studied
polytope. For
example, its extreme points are the threshold graphs. (A graph is a
threshold graph if there is a real number $S$ and for each vertex $v$ a
real vertex weight $w(v)$ such that, for any two vertices $v,u$,
$(u,v), $ there is an edge if and only if $w(u)+w(v)\ge S$. See \cite{mahadevpeled95} for much more on this.)

A self-contained proof of Theorem \ref{meanspace} is given in Section
\ref{proofmeanspace}. However, it is possible to derive it from
classical results about the mean space of exponential families (see,
e.g., \cite{brown86} or \cite{barndorffnielsen78}; in particular, see \cite{wainwrightjordan08}, Theorem~3.3).

Finally, let us describe a fast algorithm for computing the MLE if it
exists. Recall that the $L^\infty$ norm of a vector $\mathbf{x}=
(x_1,\ldots,x_n)$ is defined as
\[
|\mathbf{x}|_\infty:= \max_{1\le i\le n} |x_i|.
\]
For $1\le i\ne j\le n$ and $\mathbf{x}\in\mathbb{R}^n$, let
%
%e4 ###
\begin{equation}\label{rdef}
r_{ij}(\mathbf{x}) := \frac{1}{e^{-x_j} + e^{x_i}}.
\end{equation}
Given a realization of the random graph $G$ with degree sequence
$d_1,\ldots,d_n$, define for each $i$ the function
%
%e5 ###
\begin{equation}\label{fdef}
\varphi_i(\mathbf{x}) := \log d_i - \log\sum_{j\ne i}
r_{ij}(\mathbf{x}).
\end{equation}
Let $\varphi\dvtx \mathbb{R}^n \rightarrow\mathbb{R}^n$ be the function
whose $i$th component is
$\varphi_i$. An easy rearrangement of terms shows that the fixed points
of $\varphi$ are precisely the solutions of \eqref{mle}. The following
theorem exploits this to give an algorithm for computing the MLE in the
$\bolds{\beta}$-model.

\begin{tthm}\label{algo}
Suppose the ML equations \eqref{mle} have a solution $\hat{\bolds
{\beta}}$. Then
$\hat{\bolds{\beta}}$ is a fixed point of the function $\varphi$.
Starting from
any $\mathbf{x}_0\in\mathbb{R}^n$, define $\mathbf{x}_{k+1} =
\varphi(\mathbf{x}_k)$ for $k =
0,1,2,\ldots\,$. Then $\mathbf{x}_k$ converges to $\hat{\bolds{\beta
}}$ geometrically
fast in the~$L^\infty$ norm where the rate depends only on
$(|\hat{\bolds{\beta}}|_\infty, |\mathbf{x}_0|_\infty)$. In
particular, $\hat
{\bolds{\beta}}$
must be the unique solution of \eqref{mle}. Moreover,
\[
|\mathbf{x}_0 - \hat{\bolds{\beta}}|_\infty\le C|\mathbf{x}_0 -
\mathbf{x}_1|_\infty,
\]
where $C$
is a continuous function of the pair $(|\hat{\bolds{\beta}}|_\infty,
|\mathbf{x}_0|_\infty)$. Conversely, if the ML equations \eqref{mle}
do not
have a solution, then the sequence $\{\mathbf{x}_k\}$ must have a divergent
subsequence.
\end{tthm}

There are many other algorithms available for calculating the MLE. For
example, Holland and Leinhardt \cite{hollandleinhardt81} use an
iterative scaling algorithm and discuss the method of scoring and
weighted least squares. Hunter \cite{hunter04} develops the MM
algorithm for a similar task. Markov chain Monte Carlo algorithms and
the Robbins--Monro stochastic approximation approach are also used for
computing the MLE in exponential random graph models. See~%
\cite{kolaczyk09}, Section 6.5.2, for examples and literature. The
iterative algorithm we use is a hybrid of standard algorithms which
works well in practice and allows the strong conclusions of
Theorem~\ref{algo}. We hope that variants can be developed for related
high dimensional problems.

\begin{figure}

\includegraphics{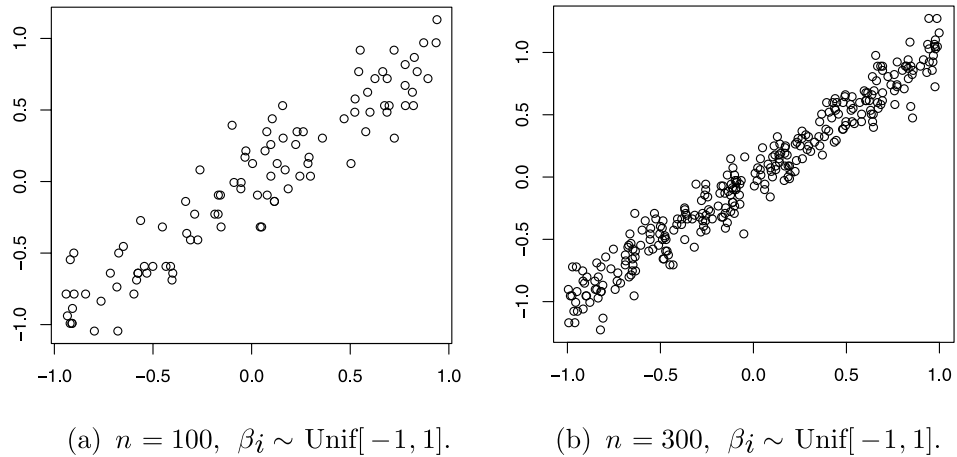}

\caption{Simulation results: plot of $\hat{\beta}_i$ vs. $\beta_i$.}\label{fig1}
\end{figure}

Let us now look at the results of some simulations. The left-hand panel
in Figure~\ref{fig1} shows the plot of $\hat{\beta}_i$ versus $\beta_i$
for a graph with $100$ vertices, where $\beta_1,\ldots,\beta_n$ were
chosen independently at uniform from the interval $[-1,1]$. The
right-hand panel is the same, except that $n$ has been increased to
$300$. The increased accuracy for larger $n$ is clearly visible.

We have also compared our results with the simulation results from the
importance sampling algorithm of Blitzstein and Diaconis
\cite{blitzsteindiaconis09} for a variety of other examples. The
results of Figure \ref{fig1} are typical. This convinces us that the
procedures developed in this paper are useful for practical problems.

%f1 ###

\textit{Comparison to the Barvinok and Hartigan work.} As
mentioned before, the present work is closely related to a recent
series of papers by Barvinok and Hartigan
\cite{barvinok08,barvinok10,barvinokHartigan09b,barvinokHartigan09c,barvinokHartigan09a,barvinokHartigan10}.
The work was initiated by Barvinok who looked at directed and bipartite
graphs in \cite{barvinok10}. In their most recent article
\cite{barvinokHartigan10} (uploaded to arXiv when our paper was near
completion), they study uniform random (undirected) graphs on $n$
vertices with a given degree sequence $\mathbf{d}= (d_1,\ldots, d_n)$ and
work with an exponential model as in Section~\ref{statsec} with~%
$\beta_i$ chosen so that the expected degree at $i$ under the
$\bolds{\beta}$-model is $d_i$. Let $G_{\mathbf{d}}$ be a uniformly
chosen random graph
with the degree sequence $\mathbf{d}$ and $G_{\bolds{\beta}}$ be a
random graph chosen
from the $\bolds{\beta}$-model. One of their main results shows that (under
hypothesis) these two graphs are close together in the following sense:
Fix a set of edges $S$ in the complete graph on $n$ vertices. Let
$X_{\mathbf{d}}$ be the number of edges of $G_{\mathbf{d}}$ in $S$.
Let $X_{\bolds{\beta}}$ be
the number of edges of $G_{\bolds{\beta}}$ in~$S$. They prove that
$X_{\mathbf{d}}/n^2$
and $X_{\bolds{\beta}}/n^2$ are each concentrated about their means (using
results from the earlier work \cite{barvinok08}) and that these means
are approximately equal. Their theorem is proved under a condition on
the degree sequences that they call ``delta tame.''

While the two sets of results (i.e., ours and those of Barvinok and
Hartigan) were proved independently and the methods of proof are quite
different in certain parts (but similar in others), the possible
connections are tantalizing. We believe that their mode of convergence
($G_{\mathbf{d}}$ and $G_{\bolds{\beta}}$ contain about the same
number of edges in a
given set) is equivalent to the graph limit convergence used here.
Perhaps this can be established using the ``cut-metric'' of Frieze and
Kannan, as expounded in \cite{borgsetal08}. We further
conjecture, based on Lemma \ref{mainlmm} in this paper, that their
delta tame condition is equivalent to our condition that the limiting
degree sequence $f$ is in the interior of~$\mathcal{F}$. If this is
so, then
Proposition~\ref{egcond} (or more accurately, Lemma~\ref{mainlmm})
gives a~necessary and sufficient condition for a degree sequence to be
delta tame, showing that essentially all degree sequences except the
ones close to the Erd\H{o}s--Gallai boundary are delta tame.

In summary, the Barvinok and Hartigan work \cite{barvinokHartigan10}
contains elegant estimates of the number of graphs with a given degree
sequence and extensions to bipartite graphs under a condition called
delta tameness; we work in the emerging language of graph limits and
prove a limit theorem under a~continuum version of the easily
verifiable Erd\H{o}s--Gallai criterion. Our work contains an efficient
algorithm for computing the maximum likelihood estimates of $\bolds
{\beta}$ for
a given degree sequence with proofs of convergence of the algorithm and
consistency of the estimates.

%The two papers may look quite different to the casual reader, but the
%underlying similarities are probably more than meets the eye.

The rest of the paper is organized as follows. In Section~\ref{proofalgo} we prove Theorem~\ref{algo}. This is followed by the
proof of Theorem~\ref{meanspace} in Section~\ref{proofmeanspace}. Both
of these theorems are required for the proof of Theorem~\ref{mainthm},
which is given in Section~\ref{proofmainthm}. Proposition~\ref{egcond}
is proved in Section~\ref{egcondproof}. Finally, the proof of Theorem~\ref{degthm}, which uses all the other theorems, is given in
Section~\ref{degthmproof}.

%s2 ###
\section{\texorpdfstring{Proof of Theorem \protect\ref{algo}}{Proof of Theorem 1.5}}\label{proofalgo}

For a matrix $A= (a_{ij})_{1\le i,j\le n}$, the $L^\infty$ operator
norm is defined as
\[
|A|_\infty:= \max_{|x|_\infty\le1} |A\mathbf{x}|_\infty.
\]
It is a simple exercise to verify that
\[
|A|_\infty= \max_{1\le i\le n} \sum_{j=1}^n |a_{ij}|.
\]
Given $\delta>0$, let us say the matrix $A$ belongs to the class
$\mathcal{L}_n(\delta)$ if $|A|_\infty\le1$ and for each $1\le i\ne
j\le n$,
\[
a_{ii} \ge\delta\quad\mbox{and}\quad a_{ij}\le-\frac{\delta}{n-1}.
\]
Lemma \ref{crucial} is our key tool.

\begin{lmm}\label{crucial}
Let $\mathcal{L}_n(\delta)$ be defined as above. If $A,B \in\mathcal
{L}_n(\delta)$,
then
\[
|AB|_\infty\le1-\frac{2(n-2)\delta^2}{n-1}.
\]
\end{lmm}

\begin{pf}
Fix $1\le i\ne k\le n$. By the definition of $\mathcal{L}_n(\delta)$,
\[
\sum_{j\notin\{ i,k\}} a_{ij}b_{jk} \ge\frac{(n-2)\delta^2}{(n-1)^2}
\quad\mbox{and}\quad a_{ii}b_{ik} + a_{ik}b_{kk}\le-\frac{2\delta^2}{n-1}.
\]
Now, if $x,y$ are two positive real numbers, then $|x-y| = |x| + |y| -
2\min\{x,y\}$. Taking $x= \sum_{j\notin\{ i,k\}} a_{ij}b_{jk}$ and
$y = -(a_{ii}b_{ik} + a_{ik}b_{kk})$, we get
\begin{eqnarray*}
\Biggl|\sum_{j=1}^n a_{ij}b_{jk}\Biggr|
&\le&
\sum_{j=1}^n |a_{ij}b_{jk}| -2\min\biggl\{\sum_{j\notin\{i,k\}} a_{ij}b_{jk},-(a_{ii}b_{ik} + a_{ik}b_{kk})\biggr\}
\\
&\le&
\sum_{j=1}^n |a_{ij}b_{jk}| - \frac{2(n-2)\delta^2}{(n-1)^2}.
\end{eqnarray*}
Combining this with the hypothesis that $|A|_\infty\le1$ and
$|B|_\infty\le1$, we get
\begin{eqnarray*}
|AB|_\infty
&=&
\max_{1\le i\le n} \sum_{k=1}^n \Biggl|\sum_{j=1}^n a_{ij}b_{jk}\Biggr|
\\
&\le&
\max_{1\le i\le n} \sum_{j=1}^n \sum_{k=1}^n |a_{ij}b_{jk}| -\frac{2(n-2)\delta^2}{n-1}
\\
&\le&
1 - \frac{2(n-2)\delta^2}{n-1}.
\end{eqnarray*}
The proof is complete.
\end{pf}

Now recall the functions $r_{ij}$ defined in \eqref{rdef}. Let
\[
q_{ij}(\mathbf{x}) := \frac{r_{ij}(\mathbf{x})}{\sum_{k\ne i}
r_{ik}(\mathbf{x})}.
\]
Note that for each $i$ and $\mathbf{x}$, $\sum_{j\ne i}
q_{ij}(\mathbf{x}) = 1$.
Again, for each $i$
\[%\label{par1}
\frac{\partial\varphi_i}{\partial x_i} = -\frac{\sum_{j\ne
i}\partial r_{ij}/\partial x_i}{\sum_{j\ne i} r_{ij}} = \sum
_{j\ne i}
\frac{e^{x_i}}{e^{-x_j} + e^{x_i}} q_{ij}
\]
and similarly for each distinct $i$ and $j$,
\[%\label{par2}
\frac{\partial\varphi_i}{\partial x_j} = -\frac
{e^{-x_j}}{e^{-x_j}+e^{x_i}} q_{ij}.
\]
Now, if $|\mathbf{x}|_\infty\le K$, then clearly
\[
\tfrac{1}{2}e^{-K} \le r_{ij}(\mathbf{x}) \le\tfrac{1}{2}e^K \qquad\mbox{for all } 1\le i\ne j\le n.
\]
Thus,
\[
\frac{e^{-2K}}{n-1} \leq q_{ij}(\mathbf{x}) = \frac{r_{ij}(\mathbf
{x})}{\sum
_{k\ne
i}r_{ik}(\mathbf{x})} \leq\frac{e^{2K}}{n-1}.
\]
It follows that for every $1\le i\ne j\le n$
%
%e6 ###
\begin{equation}\label{jbd2}
-\frac{e^{2K}}{n-1} \leq\frac{\partial\varphi_i}{\partial x_j} \le
-\frac{e^{-4K}}{2(n-1)}
\end{equation}
and also, for every $1\le i\le n$,
%
%e7 ###
\begin{equation}\label{jbd1}
\frac{1}{2}e^{-4K} \leq\frac{\partial\varphi_i}{\partial x_i} \leq e^{2K}.
\end{equation}
Now take any $\mathbf{x},\mathbf{y}\in\mathbb{R}^n$ and let $K$ be
the maximum of the
$L^\infty$ norms of $\mathbf{x}$, $\mathbf{y}$, $\varphi(\mathbf
{x})$ and $\varphi(\mathbf{y})$.
Let $J(\mathbf{x},\mathbf{y})$ be the matrix whose $(i,j){\rm{th}}$
element is
\[
J_{ij}(\mathbf{x},\mathbf{y}) = \int_0^1 \frac{\partial\varphi
_i}{\partial x_j}\bigl(t\mathbf{x}+ (1-t)\mathbf{y}
\bigr)\,dt.
\]
It is a simple calculus exercise to verify that
%
%e8 ###
\begin{equation}\label{e:JRelation}
\varphi(\mathbf{x})- \varphi(\mathbf{y}) = J(\mathbf{x},\mathbf
{y})(\mathbf{x}-\mathbf{y}).
\end{equation}
For each $i\ne j$, $\partial\varphi_i/\partial x_j$ is negative
everywhere and for each $i$, $\partial\varphi_i /\partial x_i$ is
positive everywhere. Moreover, for each $i$,
\[
\sum_{j=1}^n \biggl|\frac{\partial\varphi_i}{\partial x_j}\biggr| = \frac
{\partial\varphi_i}{\partial x_i} - \sum_{j\ne
i} \frac{\partial\varphi_i}{\partial x_j}\equiv1.
\]
It follows that for $i$ and any $\mathbf{x}, \mathbf{y}$, $\sum_{j=1}^n
|J_{ij}(\mathbf{x}, \mathbf{y})| = 1$. In particular, $|J(\mathbf
{x},\break \mathbf{y})|_\infty=
1$. From \eqref{jbd2} and \eqref{jbd1} and the fact that
$|J(\mathbf{x},\mathbf{y})|_\infty= 1$, we see that $J(\mathbf
{x},\mathbf{y})\in
\mathcal{L}_n(\delta)$ for $\delta= \frac{1}{2}e^{-4K}$. Similarly,
\begin{eqnarray*}
\varphi(\varphi(\mathbf{x})) - \varphi(\varphi(\mathbf{y}))
&=&
J(\varphi(\mathbf{x}),\varphi
(\mathbf{y})) \bigl(\varphi(\mathbf{x})-\varphi(\mathbf{y})\bigr)
\\
&=&
J(\varphi(\mathbf{x}),\varphi(\mathbf{y}))J(\mathbf{x},\mathbf
{y}) (\mathbf{x}-\mathbf{y})
\end{eqnarray*}
and $J(\varphi(\mathbf{x}),\varphi(\mathbf{y})) \in\mathcal
{L}_n(\delta)$ also. Applying Lemma
\ref{crucial}, we get
\[
|J(\varphi(\mathbf{x}), \varphi(\mathbf{y})) J(\mathbf{x},\mathbf
{y})|_\infty\le1 -
\frac{2(n-2)\delta^2}{n-1}.
\]
Thus,
%
%e9 ###
\begin{equation}\label{maineq}
|\varphi(\varphi(\mathbf{x}))-\varphi(\varphi(\mathbf
{y}))|_\infty\le\biggl(1 -
\frac{2(n-2)\delta^2}{n-1}\biggr) |\mathbf{x}- \mathbf{y}|_\infty.
\end{equation}
The quantity inside the brackets will henceforth be denoted by
$\theta(\mathbf{x},\mathbf{y})$. Note that $0\le\theta(\mathbf
{x},\mathbf{y}) < 1$ and
$\theta$ is uniformly bounded away from~$1$ on subsets of $\mathbb{R}^n
\times\mathbb{R}^n$. Moreover, since $|J(\mathbf{x}, \mathbf
{y})|_\infty=1$, we also
have the trivial but useful bound
\[
|\varphi(\mathbf{x})- \varphi(\mathbf{y})|_\infty\le|\mathbf
{x}-\mathbf{y}|_\infty.
\]
Now suppose $\varphi$ has a fixed point $\hat{\bolds{\beta}}$. If
we start with
arbitrary $\mathbf{x}_0$ and define $\mathbf{x}_{k+1} = \varphi
(\mathbf{x}_k)$ for each $k
\ge0$, then for each $k$, we have
\[
|\mathbf{x}_{k+1}-\hat{\bolds{\beta}}|_\infty= |\varphi(\mathbf
{x}_k) - \varphi(\hat
{\bolds{\beta}})|_\infty\le|\mathbf{x}_k-\hat{\bolds{\beta}}|_\infty.
\]
In particular, the sequence $\{\mathbf{x}_k\}_{k\ge0}$ remains bounded.
Therefore, by \eqref{maineq},~the\-re is a single $\theta\in[0,1)$,
depending only on $|\hat{\bolds{\beta}}|_\infty$ and $|\mathbf
{x}_0|_\infty$ in a
continuous manner, such that for all $k\ge0$, we have
%
%e10 ###
\begin{equation}\label{maineq2}
|\mathbf{x}_{k+3} - \mathbf{x}_{k+2}|_\infty\le\theta|\mathbf
{x}_{k+1}-\mathbf{x}
_k|_\infty
\end{equation}
and
\[
|\mathbf{x}_{k+2} - \hat{\bolds{\beta}}|_\infty\le\theta|\mathbf
{x}_{k}-\hat{\bolds{\beta}
}|_\infty.
\]
The second inequality shows that $\mathbf{x}_k$ converges to $\hat
{\bolds{\beta}}$
geometrically fast and the first inequality gives
\begin{eqnarray*}
|\mathbf{x}_0-\hat{\bolds{\beta}}|_\infty
&\le&
\sum_{k=0}^\infty
|\mathbf{x}_k - \mathbf{x}
_{k+1}|_\infty
\\
&\le&
\frac{1}{1-\theta} (|\mathbf{x}_0-\mathbf{x}_1|_\infty+
|\mathbf{x}_1-\mathbf{x}
_2|_\infty)
\\
&\le&
\frac{2}{1-\theta}|\mathbf{x}_0-\mathbf{x}_1|_\infty.
\end{eqnarray*}
Finally, note that if $\hat{\bolds{\beta}}$ does not exist, then the sequence
$\{\mathbf{x}_k\}$ must have a~divergent subsequence. Otherwise,
\eqref{maineq} would imply that \eqref{maineq2} must hold for all $k$
for some $\theta\in[0,1)$. This, in turn, would imply that $\mathbf{x}_k$
must converge to a limit as $k \rightarrow\infty$, which would then
be a fixed
point of $\varphi$ and, therefore, a solution of the ML equations. The
proof is complete.

Before moving to the next section we will prove a technical lemma which
will be of use in the proof of Theorem~\ref{degthm} based on the above
calculations.

\begin{lmm}\label{l:L1Recursion}
Let $\mathbf{x},\mathbf{y}\in\mathbb{R}^n$ such that $\max\{
|\mathbf{x}|_\infty,|\mathbf{y}
|_\infty\}
\leq K$. Then
\[
|\varphi(\mathbf{x}) - \varphi(\mathbf{y})|_1 \le2 e^{2K}|\mathbf
{x}- \mathbf{y}|_1,
\]
where
$|\cdot|_1$ is the usual $L^1$ norm on $\mathbb{R}^n$.
\end{lmm}

\begin{pf}
By equation~\eqref{e:JRelation},
\begin{eqnarray*}
|\varphi(\mathbf{x})- \varphi(\mathbf{y}) |_1
&=&
| J(\mathbf
{x},\mathbf{y}
)(\mathbf{x}-\mathbf{y}) |_1
\\
& =&
 \sum_{i=1}^n \Biggl| \sum_{j=1}^n (x_j-y_j)\cdot\int_0^1 \frac
{\partial\varphi_i}{\partial x_j}\bigl(t\mathbf{x}+ (1-t)\mathbf{y}\bigr)\,dt
\Biggr|
\\
&\leq&
\sum_{j=1}^n |x_j-y_j| \cdot\Biggl(\sum_{i=1}^n \sup_{t\in[0,1]} \biggl|
\frac{\partial\varphi_i}{\partial x_j}\bigl(t\mathbf{x}+ (1-t)\mathbf
{y}\bigr) \biggr|
\Biggr)
\\
&\leq&
\sum_{j=1}^n 2e^{2K}|x_j-y_j|= 2e^{2K}|\mathbf{x}- \mathbf
{y}|_1 ,
\end{eqnarray*}
where the second inequality follows from equations~\eqref{jbd2} and
\eqref{jbd1}.
\end{pf}

%s3 ###
\section{\texorpdfstring{Proof of Theorem \protect\ref{meanspace}}{Proof of Theorem 1.4}}\label{proofmeanspace}
We need the following simple technical lemma.

\begin{lmm}\label{keylmm}
Suppose $f\dvtx \mathbb{R}^n \rightarrow\mathbb{R}$ is a twice
differentiable function such that
$M := \sup_{\mathbf{x}\in\mathbb{R}^n } f(\mathbf{x}) < \infty$.
Let $\nabla f$ and
$\nabla^2 f$ denote the gradient vector and the Hessian matrix of $f$
and suppose there is a finite constant $C$ such that the $L^2$ operator
norm of $\nabla^2f$ is uniformly bounded by $C$. Then for any~$\mathbf{x}
\in
\mathbb{R}^n$,
\[
|\nabla f(\mathbf{x})|^2 \le2C \bigl(M- f(\mathbf{x})\bigr),
\]
where $|\cdot|$ denotes the Euclidean norm. In particular, there
exists a sequence $\{\mathbf{x}_{k}\}_{k\ge1}$ such that $\lim
_{k\rightarrow
\infty}
\nabla f(\mathbf{x}_k) = 0$.
\end{lmm}

\begin{pf}
Fix a point $\mathbf{x}\in\mathbb{R}^n$ and let $\mathbf{y}= \nabla
f(\mathbf{x})$. Suppose
$C$ is a uniform bound on the $L^2$ operator norm of $\nabla^2 f$. Then
for any $t\ge0$,
%
%e11 ###
\begin{equation}\label{gradbd}
|\nabla f(\mathbf{x}+ t\mathbf{y}) - \nabla f(\mathbf{x}) | \le
Ct|\mathbf{y}|.
\end{equation}
Now let $g(t) = f(\mathbf{x}+ t\mathbf{y})$. Then for all $t$,
\[
g(t) - g(0) \le M - f(\mathbf{x}).
\]
Again, note that
\begin{eqnarray*}
g'(t)
&=&
 \langle\mathbf{y}, \nabla f(\mathbf{x}+ t\mathbf
{y})\rangle
\\
&=&
 \langle\mathbf{y}, \nabla f(\mathbf{x}+ t\mathbf{y}) - \nabla
f(\mathbf{x})\rangle+
\langle\mathbf{y}, \nabla f(\mathbf{x})\rangle
\\
&\ge&
- Ct|\mathbf{y}|^2 + |\mathbf{y}|^2.
\end{eqnarray*}
[The last step follows by \eqref{gradbd} and Cauchy--Schwarz.] Thus,
for any $t\ge0$,
\[
M - f(\mathbf{x}) \ge g(t) - g(0) = \int_0^t g'(s)\,ds \ge|\mathbf
{y}|^2 \int_0^t
(1-Cs)\,ds.
\]
Taking $t = 1/C$ gives the desired result.
\end{pf}

\begin{pf*}{Proof of Theorem \ref{meanspace}}
Let $g = (g_1,\ldots,g_n)\dvtx \mathbb{R}^n \rightarrow\mathbb{R}^n$ be
the function defined as
\[
g_i(\mathbf{x}) = \sum_{j\ne i} \frac
{e^{x_i+x_j}}{1+e^{x_i+x_j}},\qquad i=1,\ldots,
n.
\]
Then $\mathcal{R}$ is the range of $g$. This is because the expected
degree of
vertex $i$ of a random graph following the law $\mathbb{P}_\mathbf
{x}$ is
$g_i(\mathbf{x})$. In particular, the vector $g(\mathbf{x})$ is a
weighted average
of degree sequences and hence,
\[
\operatorname{conv}({\mathcal{D}}) \supseteq\overline{\mathcal{R}}.
\]
Now, for every $\mathbf{y}\in\mathbb{R}^n$, let $f_\mathbf
{y}\dvtx \mathbb{R}^n \rightarrow\mathbb{R}$ be the
function
\[
f_\mathbf{y}(\mathbf{x}) = \sum_{i=1}^n x_i y_i - \log\sum_{1\le
i<j\le n}
(1+e^{x_i + x_j}).
\]
Note that under $\mathbb{P}_\mathbf{x}$, the probability of obtaining
a given graph
with degree sequence $d = (d_1,\ldots,d_n)$ is exactly
\[
\frac{e^{\sum_ix_i d_i}}{\prod_{i<j} (1+e^{x_i +x_j})}.
\]
Thus, the above quantity must be bounded by $1$ and hence, taking logs,
we get $f_d(\mathbf{x}) \le0$. Since $f_\mathbf{y}(\mathbf{x})$
depends linearly on
$\mathbf{y}$, this implies that
\[
f_\mathbf{y}(\mathbf{x}) \le0\qquad \mbox{for all } \mathbf{y}\in
\operatorname{conv}(\mathcal{D}), \mathbf{x}\in\mathbb{R}^n.
\]
Now fix $\mathbf{y}\in\operatorname{conv}(\mathcal{D})$. Then $f_\mathbf
{y}(\mathbf{x})\le0$ for all
$\mathbf{x}\in
\mathbb{R}^n$. Moreover,
%[-1,0],
it is easy to show that $\nabla^2 f$ is uniformly bounded. Hence, it
follows from Lemma~\ref{gradbd} that there exists a sequence
$\{\mathbf{x}_k\}_{k\ge1}$ such that $\lim_{k\rightarrow\infty}
\nabla
f_\mathbf{y}(\mathbf{x}_k) =0$. But
\[
\nabla f_\mathbf{y}(\mathbf{x}) = \mathbf{y}- g(\mathbf{x}).
\]
Thus, $\mathbf{y}= \lim_{k\rightarrow\infty} g(\mathbf{x}_k)$.
This shows that
\[
\operatorname{conv}(\mathcal{D}) \subseteq\overline{\mathcal{R}}
\]
and hence, completes the proof of the claim that $\operatorname{conv}(\mathcal{D}) =
\overline{\mathcal{R}}$.
\end{pf*}
%
%Finally, to see that $\mr$ is not closed, simply consider the complete
%graph on $n$ vertices, which gives the degree sequence $(n-1,
%range of $g$.

%s4 ###
\section{\texorpdfstring{Proof of Theorem \protect\ref{mainthm} (Existence and consistency of the MLE)}%
{Proof of Theorem 1.3 (Existence and consistency of the MLE)}}\label{proofmainthm}

This section uses the notation of Section \ref{intro} without explicit
reference. The proof consists of two lemmas. The first lemma gives a
condition for the ``tightness'' of the MLE. This result is closely
related to the Erd\H{o}s--Gallai characterization of degree sequences.
The second lemma shows that conditions needed for the first lemma are
satisfied with high probability. An addenda at the end of the section
contains some results about existence of the MLE and the closely
related topic of conjugate Bayesian analysis.

In this section we will repeatedly encounter statements like ``$c$ is a
positive constant depending only on $a, b,\ldots\,$.'' Such a statement
should be interpreted as ``$c$ can be expressed as a function of $a,
b,\ldots$ that is bounded away from~$0$ and $\infty$ on compact subsets
of the domain of $(a, b,\ldots).$'' Sometimes, $c$ will be expressed as
$C(a,b,\ldots)$.

\begin{lmm}\label{mainlmm}
Let $(d_1,\ldots,d_n)$ be a point in the set $\overline{\mathcal
{R}}$ of
Theorem \ref{meanspace}. Suppose there exist $c_1, c_2\in(0,1)$ such
that $c_2(n-1) \le d_i \le c_1 (n-1)$ for all~$i$. Suppose $c_3$ is a
positive constant such that %Then there exists $b\in(0,1)$ depending
%only on $c_1,c_2$, such that if the quantity
%
\[
\frac{1}{n^2}\inf_{B\subseteq\{1,\ldots,n\}, |B|\ge
c_2^2n}\biggl\{\sum_{j\notin B} \min\{d_j, |B|\} + |B|(|B|-1) - \sum
_{i\in
B} d_i \biggr\} \ge c_3.
\]
Then a solution $\hat{\bolds{\beta}}$ of \eqref{mle} exists and satisfies
$|\hat{\bolds{\beta}}|_\infty\le c_4$, where $c_4$ is a constant
that depends
only on $c_1,c_2,c_3$.
\end{lmm}

\begin{pf}
In this proof, $C(c_1,c_2,c_3)$ denotes positive constants that depend
only on $c_1, c_2, c_3$, in the sense defined above. The argument
repeatedly uses the monotonicity of $e^{x+y}/(1+e^{x+y})$ in $x$ for
each $y$.

Assume first that $\hat{\bolds{\beta}}$ exists in the sense that
there exists
$\hat{\bolds{\beta}}\in\mathbb{R}^n$ such that~\eqref{mle} is
satisfied. Let $c_1,
c_2, c_3$ be as in the statement of the lemma. It is proved below that
$|\hat{\bolds{\beta}}|_\infty$ is bounded above by $C(c_1,c_2,c_3)$.

Let $d_{\max}:= \max_i d_i$ and $d_{\min} := \min_i d_i$. Similarly,
let $\hat{\beta}_{\max} := \max_i \hat{\beta}_i$ and $\hat{\beta
}_{\min} := \min_i \hat{\beta}
_i$. The
first step is to prove that $\hat{\beta}_{\max}\le C(c_1,c_2,c_3)$. If
\mbox{$\hat{\beta}_{\max}\le0$}, there is nothing to prove. So assume that
$\hat{\beta}_{\max} > 0$. Let
\[
m := \bigl|\bigl\{i\dvtx \hat{\beta}_i > -{\textstyle\frac{1}{2}} \hat{\beta
}_{\max}\bigr\}\bigr|.
\]
Clearly, by the assumption that $\hat{\beta}_{\max}>0$, it is
guaranteed that
$m \ge1$. Let~$i^*$ be an index that maximizes $\hat{\beta}_i$. Then by
\eqref{mle}, we see that
\[
d_{\max}\ge d_{i^*} > (m-1)
\frac{e^{(1/2)\hat{\beta}_{\max}}}{1+e^{(1/2)\hat
{\beta}_{\max}}}.
\]
This implies
\[
n - m > n -1 - d_{\max}\bigl(1+e^{-(1/2)\hat{\beta}_{\max}}\bigr)
\ge
n - 1 -c_1 (n-1) \bigl(1+e^{-(1/2)\hat{\beta}_{\max}}\bigr).
\]
In particular, this shows that if $\hat{\beta}_{\max} > C(c_1)$ then
$m < n$
and hence, there exists $i$ such that $\hat{\beta}_i \le
-\frac{1}{2}\hat{\beta}_{\max}$. Suppose this is true and fix any
such $i$.
(In particular note that $\hat{\beta}_i < 0$.) Let
\[
m_i := \bigl|\bigl\{j \dvtx j\ne i, \hat{\beta}_j < -{\textstyle\frac{1}{2}}\hat
{\beta}_i\bigr\}\bigr|.
\]
Then by \eqref{mle},
\[
d_{\min} \le d_i < m_i
\frac{e^{(1/2)\hat{\beta}_i}}{1+e^{(1/2)\hat{\beta
}_i}} + n-1 - m_i,
\]
which gives
\[
m_i < (n-1-d_{\min}) \bigl(1+e^{(1/2)\hat{\beta}_i}\bigr)
\le
(n-1)(1-c_2) \bigl(1+e^{-(1/4)\hat{\beta}_{\max}}\bigr).
\]
Note that there are at least $n-m_i$ indices $j$ such that $\hat{\beta
}_j \ge
-\frac{1}{2}\hat{\beta}_i \ge\frac{1}{4}\hat{\beta}_{\max}$.
The last display implies that if $\hat{\beta}_{\max} > C(c_1, c_2)$, then
there exists $i$ such that $n-m_i \ge bn$, where %$b\in(0,1)$ is a
%constant that depends only on $c_1, c_2$.
%
\[
b := c_2^2.
\]
Consequently, if $\hat{\beta}_{\max} > C(c_1, c_2)$, there is a set
$A\subseteq\{1,\ldots,n\}$ of size at least $bn$ such that $\hat
{\beta}_j
\ge\frac{1}{4}\hat{\beta}_{\max}$ for all $j\in A$, where $b=
c_2^2$. %$b
Henceforth, assume that~$\hat{\beta}_{\max}$ is so large that such a set
exists. Let
\[
h := \sqrt{\hat{\beta}_{\max}}.
\]
For each integer $r$ between $0$ and $\frac{1}{16}h - 1$, let
\[
D_r := \biggl\{i\dvtx -\frac{\hat{\beta}_{\max}}{8} + rh \le\hat
{\beta}_i <
-\frac{\hat{\beta}_{\max}}{8} + (r+1)h\biggr\}.
\]
Since $D_0,D_1,\ldots$ are disjoint, there exists $r$ such that
\[
|D_r|\le\frac{n}{(1/16)h-1},
\]
provided $h > 16$. By assumption, $\hat{\beta}_{\max} > C(c_1,c_2)$.
Since we
are free to choose $C(c_1,c_2)$ as large as we like, it can be assumed
without loss of generality that $h > 16$.

Fix such an $r$ between $0$ and $\frac{1}{16}h - 1$. Let
\[
B := \biggl\{ i\dvtx \hat{\beta}_i \ge\frac{\hat{\beta}_{\max
}}{8} - \biggl(r +
\frac{1}{2}\biggr) h\biggr\}.
\]
Clearly, the set $B$ contains the previously defined set $A$ and hence,
%
%e12 ###
\begin{equation}\label{bsize}
|B|\ge b n.
\end{equation}
Now, for each $i\ne j$, define
\[
%p_{ij} := \frac{e^{\beta_i + \beta_j}}{1+e^{\beta_i +\beta_j}}, \ \ \
\hat{p}_{ij} := \frac{e^{\hat{\beta}_i + \hat{\beta
}_j}}{1+e^{\hat{\beta}_i +\hat{\beta}_j}}.
\]
For each $i$, let
\[
d_i^B := \sum_{j\in B\backslash\{i\}} \hat{p}_{ij}.
\]
Since $\hat{\beta}_i \ge\frac{\hat{\beta}_{\max}}{16}$ for each
$i\in B$, it follows
that
%
%e13 ###
\begin{eqnarray}\label{ineq1}
|B|(|B|-1) - \sum_{i\in B} d_i^B
&=&
|B|(|B|-1) - \sum_{i,j\in B,i\ne j} \hat{p}_{ij}\nonumber
\\
&=&
\sum_{i,j\in B, i\ne j} (1-\hat{p}_{ij})
\\
&\le&
\frac{|B|(|B|-1) }{1+e^{(1/8)\hat{\beta}_{\max}}}.\nonumber
\end{eqnarray}
The above inequality is the first step of a two-step argument. For the
second step, take any $j\notin B$. Consider three cases. First,
suppose $\hat{\beta}_j \ge-\frac{\hat{\beta}_{\max}}{8} + (r+1)
h$. Then for each
$i\in B$, $\hat{\beta}_i + \hat{\beta}_j \ge\frac{h}{2}$ and, therefore,
\[
\min\{d_j, |B|\} - d_j^B \le|B|- \sum_{i\in B} \hat{p}_{ij}
= \sum_{i\in B} (1-\hat{p}_{ij})\le\frac{|B|}{1+e^{h/2}}.
\]
Next, suppose $\hat{\beta}_j \le-\frac{\hat{\beta}_{\max}}{8} +
rh$. Then for any
$i\notin B$, $\hat{\beta}_i + \hat{\beta}_j \le-\frac{h}{2}$. Thus,
\[
\min\{d_j, |B|\} - d_j^B \le d_j - d_j^B
= \sum_{i\notin B, i\ne j} \hat{p}_{ij}\le n e^{-h/2}.
\]
Finally, the third case covers all $j\notin B$ that do not fall in
either of the previous two cases. This is a subset of the set of all
$j$ comprising the set $D_r$. Combining the three cases gives
%
%e14 ###
\begin{equation}\label{ineq2}
\sum_{j\notin B}(\min\{d_j, |B|\} - d_j^B) \le
\frac{n^2}{1+e^{h/2}} + n^2 e^{-h/2} +
\frac{16n^2}{h-16}.
\end{equation}
But
\[
\sum_{j\notin B} d_j^B = \sum_{i\in B, j\notin B} \hat{p}_{ij} =
\sum_{i\in B} (d_i - d_i^B).
\]
Thus, adding \eqref{ineq1} and \eqref{ineq2},
%
%e16 ###
%e15 ###
\begin{eqnarray}\label{gineq}
&&
\sum_{j\notin B} \min\{d_j, |B|\} + |B|(|B|-1) - \sum_{i\in B}
d_i\nonumber
\\[-8pt]\\[-8pt]
&&\qquad\le
\frac{n^2 }{1+e^{(1/8)\hat{\beta}_{\max}}} +\frac{n^2}{1+e^{h/2}} + n^2 e^{-h/2} +\frac{16n^2}{h-16}.\nonumber
\end{eqnarray}
The left-hand side of the above inequality is bounded below by $c_3
n^2$, by the definition of $c_3$ in the statement of the theorem. The
coefficient of $n^2$ on the right-hand side tends to zero as
$\hat{\beta}_{\max}\rightarrow\infty$. This shows that $\hat
{\beta}_{\max}\le C(c_1, c_2,
c_3)$, where the bound is finite since $c_3 >0$. Next, note that for
any~$i$,
\[
d_i \le\frac{ne^{\hat{\beta}_i + \hat{\beta}_{\max}}}{1+e^{\hat
{\beta}_i + \hat{\beta}_{\max}}}
\]
and, therefore, if $i^{**}$ is a vertex that minimizes $\hat{\beta
}_i$, then
\[
d_{\min} \le d_{i^{**}} \le\frac{ne^{\hat{\beta}_{\min} +
\hat{\beta}_{\max}}}{1+e^{\hat{\beta}_{\min} + \hat{\beta
}_{\max}}}.
\]
Combined with the upper bound on $\hat{\beta}_{\max}$ and the lower
bound on
$d_{\min}$, this shows that $\hat{\beta}_{\min} \ge-C(c_1,c_2,c_3)$.

To complete the proof of the lemma, it must be proved that $\hat
{\bolds{\beta}}$
exists. Since $(d_1,\ldots, d_n)\in\overline{\mathcal{R}}$, by Theorem
\ref{meanspace} there is a sequence of points $\{\mathbf{x}_k\}_{k\ge
0}$ in
$\mathbb{R}^n$ that converge to $(d_1,\ldots,d_n)$ for which solutions
to~\eqref{mle} exist. Let $\{\hat{\bolds{\beta}}_k\}_{k\ge0}$
denote a sequence
of solutions. The steps above prove that $|\hat{\bolds{\beta
}}_k|_\infty\le C$
for all large enough $k$ where $C$ is some constant depending only on
$c_1,c_2,c_3$. Therefore, the sequence $\{\hat{\bolds{\beta}}_k\}
_{k\ge0}$ must
have a limit point. This limit point is clearly a solution to
\eqref{mle} for the original sequence $d_1,\ldots, d_n$.
\end{pf}

The next lemma shows that the degree sequence in a typical realization
of our random graph satisfies the conditions of Lemma \ref{mainlmm}.

\begin{lmm}\label{condlmm}
Let $G$ be drawn from the probability measure $\mathbb{P}_{\bolds
{\beta}}$ and let
$d_1,\ldots, d_n$ be the degree sequence of $G$. Let $L := \max_{1\le
i\le n} |\beta_i|$ and let $c\in(0,1)$ be any constant. Then there are
constants $C > 0$ and $ c_1,c_2\in(0,1)$ depending only on $L$ and a
constant $c_3\in(0,1)$ depending only on $L$ and $c$ such that if $n >
C$, then with probability at least $1- 2n^{-2}$, $c_2(n-1)\le d_i \le
c_1(n-1)$ for all $i$ and
\begin{eqnarray*}
&&
\frac{1}{n^2}\inf_{B\subseteq\{1,\ldots,n\}, |B|\ge c n}\biggl\{
\sum_{j\notin B} \min\{d_j, |B|\} + |B|(|B|-1) - \sum_{i\in B} d_i
\biggr\}
\\
&&\qquad\ge
c_3 -\sqrt{\frac{6\log n}{n}}.
\end{eqnarray*}
\end{lmm}

\begin{pf}
Let
\[
\bar{d}_i := \sum_{j\ne i}\frac{e^{\beta_i + \beta_j}}{1+e^{\beta
_i +
\beta_j}}.
\]
%
%Let $\db$ and $\bar{\db}$ denote the vectors $(d_1,\ldots,d_n)$ and $(
Note that for each $i$, $d_i$ is a sum of independent indicator random
variables and $\mathbb{E}(d_i)=n \bar{d}_i$. Therefore, by Hoeffding's
inequality \cite{hoeffding63},
\[
\mathbb{P}(|d_i - \bar{d}_i| > x) \le2e^{-x^2/2n}.
\]
Thus, if we let $E$ be the event
\[
\Bigl\{\max_i |d_i - \bar{d}_i| > \sqrt{6n\log n}\Bigl\},
\]
then by a union bound,
\[
\mathbb{P}(E) \le\frac{2}{n^2}.
\]
Now, clearly, there are constants $c_1'<1 $ and $c_2'>0$ depending only
on $L$ such that $c_2'(n-1)\le\bar{d}_i \le c_1'(n-1)$ for all $i$.
Therefore, under $E^c$, if $n$ is sufficiently large (depending on
$L$), we get constants $c_1, c_2$ depending only on $L$ such that
$c_2(n-1)\le d_i \le c_1(n-1)$ for all $i$.

Next, define
\[
g(d_1,\ldots,d_n, B):= \sum_{j\notin B} \min\{d_j, |B|\} + |B|(|B|-1)
- \sum_{i\in B} d_i.
\]
Note that
\[
|g(d_1,\ldots,d_n, B) - g(\bar{d}_1,\ldots,\bar{d}_n, B)| \le
\sum_{i=1}^n |d_i - \bar{d}_i| \le n \max_i\!|d_i - \bar{d}_i|.
\]
Moreover, following the notation introduced in the proof of Lemma
\ref{mainlmm}, we have
\begin{eqnarray*}
&&
g(\bar{d}_1,\ldots,\bar{d}_n, B)
\\
&&\qquad= \sum_{j\notin B} (\min\{\bar{d}_j, |B|\} - \bar{d}_j^B) +
|B|(|B|-1) - \sum_{i\in B} \bar{d}_i^B
\\
&&\qquad\ge
|B|(|B|-1) - \sum_{i\in B} \bar{d}_i^B
\\
&&\qquad=
\sum_{i,j\in B, i\ne j} (1-p_{ij})\ge c_4 |B|(|B|-1),
\end{eqnarray*}
where $c_4\in(0,1)$ is a constant depending only on $L$. Thus, under
$E^c$, $n>C$ and $|B|\ge cn$ we have
\[
g(d_1,\ldots, d_n, B) \ge c_3n^2-n^{3/2}\sqrt{6\log n},
\]
where $c_3\in(0,1)$ is a constant depending only on $L$ and $c$. The
proof is complete.
\end{pf}

\begin{pf*}{Proof of Theorem \ref{mainthm}}
Let $E$ be the event defined in the proof of Lemma~\ref{condlmm}. Let
$C, c_1,c_2$ be as in Lemma~\ref{condlmm}. By lemmas, if $E^c$ happens and $n >C$, then a solution
$\hat{\bolds{\beta}}$ of \eqref{mle} exists and satisfies
$|\hat{\bolds{\beta}}|_\infty\le C(L)$, where~$C(L)$ generically
denotes a
constant that depends only on $L$. This proves the existence of the
MLE. The uniqueness follows from Theorem \ref{algo}.

The proof of the error bound uses Theorem \ref{algo}. Let $\mathbf
{x}_0 =
\bolds{\beta}$ and define $\{\mathbf{x}_k\}_{k\ge1}$ as in Theorem
\ref{algo}. A
simple computation shows that the $i$th component of $\mathbf{x}_0 -
\mathbf{x}_1$
is simply $\log(\bar{d}_i/d_i)$. Under $E^c$ and $n >C$, this is
bounded by $C(L)\sqrt{n^{-1}\log n}$. The error bound now follows
directly from Theorem~\ref{algo}.

Finally, to remove the condition $n >C$, we simply increase $C(L)$ in
Theorem~\ref{mainthm} so that $1-C(L)n^{-2} < 0$ for $n \le C$. The
proof of Theorem~\ref{mainthm} is complete.
\end{pf*}

\textit{Addenda.} (A) \textit{Practical remarks on the MLE.} Theorem
\ref{mainthm} shows that with high probability, under the $\mathbb
{P}_{\bolds{\beta}}$
measure, for large $n$ the MLE exists and is unique. In applications, a
graph is given and Theorem \ref{mainthm} may be used to test the
$\mathbb{P}_{\bolds{\beta}}$ model. The MLE may fail to exist because
the maximum is
taken on at $\beta_i = \pm\infty$ for one or more values of $i$. For
example, with $n= 2$ vertices, an observed graph will either have zero
edges or one edge. In the first case, the likelihood is
$1/(1+e^{\beta_1+\beta_2})$, maximized at $\beta_1 = \beta_2 =
-\infty$. In the second case the likelihood is $e^{\beta_1
+\beta_2}/(1+e^{\beta_1+\beta_2})$, maximized at $\beta_1 =
\beta_2=\infty$. Here, the MLE fails to exist with probability one.

Similar considerations hold when the observed graph has any isolated
vertices and for a star graph. We conjecture: Let $G$ be a graph on $n$
vertices. The MLE for the $\bolds{\beta}$-model exists if and only if
the degree
sequence lies in the interior of the convex polytope $\operatorname{conv}(\mathcal{D})$
defined in Theorem \ref{meanspace}.

In cases where the MLE does not exist, it is customary to add a small
amount to each degree (see the discussion in \cite{bishopetal75}). This is often done in a~convenient and principled
way by using a Bayesian argument.

(B) \textit{Conjugate prior analysis for the $\bolds{\beta}$-model.}
Background on conjugate priors for exponential families is in
\cite{diaconisylvisaker79} and \cite{smithetal95,smithetal96}. The
$\bolds{\beta}$-model
\[
\mathbb{P}_{\bolds{\beta}}(G) = Z(\bolds{\beta})^{-1}e^{\sum
_{i=1}^n d_i(G) \beta_i},\qquad \bolds{\beta}\in
\mathbb{R}^n, Z(\bolds{\beta}) = \prod_{1\le i<j\le n} (1+e^{\beta
_i + \beta_j})
\]
has sufficient statistic $\mathbf{d}= (d_1,\ldots,d_n)$. Here
$\mathbf{d}$ takes
values in $\mathcal{D}$, the set of degree sequences for graphs on $n$
vertices. Thus, $\mathbb{P}_{\bolds{\beta}}$ induces a natural
exponential family on
$\mathcal{D}$ with a base measure $\mu$ that does not depend on
$\bolds{\beta}$.
Following notation in \cite{diaconisylvisaker79}, write
%
%e17 ###
\[
\mathbb{P}_{\bolds{\beta}}(\mathbf{d}) = \mu(\mathbf{d}) e^{\bolds
{\beta}\cdot\mathbf{d}- m(\bolds{\beta})}
\qquad\mbox{with }
m(\bolds{\beta}) = \log Z(\bolds{\beta}) = \sum_{1\le i<j\le n}
\log(1+e^{\beta_i +
\beta_j}).
\]
Following \cite{diaconisylvisaker79}, for $\mathbf{d}_0$ in the
interior of
$\operatorname{conv}(\mathcal{D})$ and $n_0> 0$, define the conjugate prior
of $\mathbb{R}^n$ by
\[
\pi_{n_0,\mathbf{d}_0}(\bolds{\beta}) = Z(n_0, \mathbf
{d}_0)^{-1}e^{n_0\mathbf{d}_0 \cdot\bolds{\beta}-
n_0m(\bolds{\beta})}.
\]
Here $Z(n_0, \mathbf{d}_0)$ is the normalizing constant, shown to be positive
and finite in \cite{diaconisylvisaker79}. By the theory in
\cite{diaconisylvisaker79}, $\nabla m(\bolds{\beta}) = \mathbb
{E}_{\bolds{\beta}}(\mathbf{d})$ and
\[
\mathbb{E}_{\pi_{n_0,\mathbf{d}_0}}(\nabla m (\bolds{\beta})) =
\mathbb{E}_{\pi_{n_0,\mathbf{d}_0}}(\mathbb{E}_{\bolds{\beta}}(\mathbf{d})) = \mathbf{d}_0.
\]
This identity characterizes the prior $\pi_{n_0,\mathbf{d}_0}$. The posterior,
given an observed degree sequence $\mathbf{d}(G)$, is
\[
\pi_{n_0+1, (\mathbf{d}(G)+n_0\mathbf{d}_0)/(n_0+1)}.
\]
Clearly, the mode of the posterior can be found by using the iteration
of Theorem \ref{algo}. The proof of Theorem \ref{algo} shows that the
mode exists uniquely for any observed $\mathbf{d}(G)$. The posterior
mean must
be found using standard Markov chain Monte Carlo techniques.

A natural way to obtain feasible prior mean parameters [i.e., values
of~$\mathbf{d}_0$ that lie within the interior of $\operatorname{conv}(\mathcal{D})$] is to consider a
model of random graphs that puts positive mass on every possible graph
on $n$ vertices and take its expected degree sequence. For example, the
Erd\H{o}s--R\'enyi graph $G(n,p)$, for $0<p<1$, is one such model. Its
expected degree sequence is $(c,c,\ldots, c)$ where $c= (n-1)p$. Thus,
$(c,c,\ldots,c)$ is a feasible mean parameter for every $c\in(0,n-1)$.
Similarly, the expected degree sequence in any of the standard models
of power law graphs is a feasible value of $\mathbf{d}_0$ that has
power law
behavior.

%s5 ###
\section{\texorpdfstring{Proof of Proposition \protect\ref{egcond} (characterization of the interior)}%
{Proof of Proposition 1.2 (characterization of the interior)}}\label{egcondproof}
\mbox{}
\begin{pf*}{Proof of Proposition \protect\ref{egcond}}
Let us begin by restating the Erd\H{o}s--Gallai criterion from
Section~\ref{intro}. Suppose $d_1 \ge d_2 \ge\cdots\ge d_n$ are
nonnegative integers. The Erd\H{o}s--Gallai criterion says that
$d_1,\ldots,d_n$ can be the degree sequence of a simple graph on $n$
vertices if and only if $\sum_{i=1}^n d_i$ is even and for each $1\le
k\le n$,
\[
\sum_{i=1}^k d_i \le k(k-1) + \sum_{i=k+1}^n \min\{d_i, k\}.
\]
Now take any function $f\in{D'[0,1]}$ and let
\[
G_f(x) := \int_x^1 \min\{f(y), x\}\,dy + x^2 - \int_0^x f(y)\,dy.
\]
Clearly, $G_f(x)$ is continuous as a function of $x$. If $f\in\mathcal{F}$,
the E--G criterion clearly shows that $G_f$ must be a nonnegative
function. We claim that this implies that if $f$ belongs to the
interior of $\mathcal{F}$, then $G_f(x)$ must be strictly positive for every
$x\in(0,1]$. Otherwise, there exists $x\in(0,1]$ such that $G_f(x) =
0$. If we show that there exists a sequence $f_n \rightarrow f$ in the modified
$L^1$ topology such that $G_{f_n}(x) < 0$ for each $n$, then we get a
contradiction which proves the claim. This is quite easily done by
producing $f_n$ that is strictly bigger than $f$ in $[0,x)$ and equal
to $f$ elsewhere, all the while maintaining left-continuity.

Similarly, it is clear that any $f\in\mathcal{F}$ must take values in $[0,1]$.
If $f$ attains $0$ or $1$, then we can produce a sequence $f_n
\rightarrow f$
whose ranges are not contained in $[0,1]$ and, therefore, $f$ cannot
belong to the interior of $\mathcal{F}$.

Thus, we have proved that if $f$ belongs to the interior of $\mathcal
{F}$, then
$f$ must satisfy the two conditions of Proposition \ref{egcond}. Let us
now prove the converse. Suppose $f\in{D'[0,1]}$ such that $0 < c_1 <
f(x) < c_2 < 1$ for all $x\in[0,1]$ and $G_f(x) > 0$ for all $x\in
(0,1]$. We have to show that any function that is sufficiently close to
$f$ in the modified $L^1$ norm must belong to $\mathcal{F}$.

To do that let us first prove that $f\in\mathcal{F}$. Take any $n$.
Let $d_i^n
= \lfloor n f(i/n)\rfloor$, $i=2,\ldots, n$, and $d_1^n = \lfloor n
f(0)\rfloor$. Since $f$ is nonincreasing, we have $d_1^n \ge d_2^n \ge
\cdots\ge d_n^n$. Increase some of the $d_i^n$'s by $1$, if necessary,
so that $\sum d_i^n$ is even (and monotonicity is maintained). With
this construction, it is clear that
\[
\biggl|\frac{d_1^n}{n} - f(0)\biggr| + \biggl| \frac{d_n^n}{n} - f(1)\biggr| + \frac{1}{n}
\sum_{i=1}^n\biggl|\frac{d_i^n}{n} - f\biggl(\frac{i}{n}\biggr)\biggr|\le\frac{4}{n}.
\]
Thus, if $\mathbf{d}^n$ denotes the vector $(d_1^n,\ldots, d_n^n)$, then
$\mathbf{d}^n$ converges to the scaling limit $f$. We need to show
that for
all large enough $n$, $\mathbf{d}^n$ is a valid degree sequence.

Since $f$ is bounded and nonincreasing,\vspace*{-2pt}
\[
\lim_{n\to\infty} \int_0^1 |f(x) - f(\lceil n x \rceil/n)|\,dx = 0\vspace*{-2pt}
\]
and so uniformly in $1\leq k \leq n$,
\[
\biggl|\frac{\sum_{i=k+1}^n \min\{d_i^n, k\} + k(k-1) - \sum_{i=1}^k
d_i^n}{n^2} - G_f(k/n)\biggr| \le\varepsilon(n),
\]
where $\varepsilon(n) \rightarrow0$ as $n \rightarrow\infty$. Thus,
there exists a sequence
of integers $\{k_0(n)\}$, where $k_0(n)/n \rightarrow0$ as $n
\rightarrow\infty$,
such that whenever $k \ge k_0(n)$, we have\vspace*{-2pt}
\[
\sum_{i=k+1}^n \min\{d_i^n, k\} + k(k-1) - \sum_{i=1}^k d_i^n > 0.\vspace*{-2pt}
\]
Again, there exists $c_1' <1$ and $c_2' > 0$ such that if $n$ is
sufficiently large, we have $c_2' \le d_i^n/n \le c_1'$ for all $i$.
Suppose $n$ is so large that $k_0(n)/n < c_2'$ and $(1-c_1')n- k_0(n) >
0$. Then, if $k \le k_0(n)$, we have\vspace*{-2pt}
\begin{eqnarray*}
&&\sum_{i=k+1}^n \min\{d_i^n, k\} + k(k-1) - \sum_{i=1}^k d_i^n
\\[-2pt]
&&\qquad\ge\sum_{i=k+1}^n \min\{c_2'n, k\} + k(k-1) - \sum_{i=1}^k nc_1'
\\[-2pt]
&&\qquad= (n-k)k + k(k-1) - c_1'nk
\\[-2pt]
&&\qquad= \bigl((1-c_1') n - k\bigr) k + k(k-1) >0.
\end{eqnarray*}
Thus, for $n$ so large, we have that for all $1\le k\le n$,\vspace*{-2pt}
\[
\sum_{i=k+1}^n \min\{d_i^n, k\} + k(k-1) - \sum_{i=1}^k d_i^n >0.\vspace*{-2pt}
\]
By the Erd\H{o}s--Gallai criterion, this shows that $(d_1^n,\ldots,
d_n^n)$ is a valid degree sequence.

Thus, we have shown that any $f$ that satisfies the two conditions of
Proposition~\ref{egcond} must belong to $\mathcal{F}$. Now we only
have to show
that if $f$ satisfies the two criteria, then any $h$ sufficiently close
to $f$ in the modified $L^1$ norm must also satisfy them.

Note that $G_f$ is a continuous function that is positive in $(0,1]$.
Moreover, for all $0\leq x \leq1$,\vspace*{-2pt}
\[
|G_f(x)-G_{f'}(x)|\leq\|f-f'\|_{1'}\vspace*{-2pt}
\]
so if $f_n \rightarrow f$ in the modified $L^1$ norm, then
$G_{f_n}\rightarrow G_f$ in
the supnorm. Thus, for any $\varepsilon> 0$ there exists $\delta>0$ such that
whenever $\|h-f\|_{1'}< \delta$, we have $G_h(x) > 0$ for all $x \in
[\varepsilon,1]$. We also have that $c_1-\delta\leq h(x) \leq c_2
+\delta$ for
all $0\leq x \leq1$. Choosing $\delta,\varepsilon>0$ small as
necessary, we
can ensure that $c_1 -\delta> \varepsilon$ and $1-\varepsilon-\delta- c_2
> 0$. Fix
such $\varepsilon$, $\delta$ and $h$. Then, for $x\in(0,\varepsilon)$, we have\vspace*{-2pt}
\begin{eqnarray*}
G_h(x)
&\ge&
\int_x^1 \min\{c_1 - \delta, x\}\,dy + x^2 - \int_0^x(c_2 + \delta)\,dy
\\[-2pt]
&=&
(1-x)x + x^2 - (c_2 + \delta) x
\\[-2pt]
&=&
(1-\varepsilon-\delta- c_2) x + x^2 >0.
\end{eqnarray*}
But we also have $G_h(x) > 0$ for $x\in[\varepsilon,1]$ by the choice of
$\delta$. Thus, we have proved that there exists $\delta>0$ such that
whenever $\|h-f\|_{1'}< \delta$, we have $G_h(x) >0$ for all $x\in
(0,1]$. Choosing $\delta$ sufficiently small, we can ensure that the
range of $h$ does not contain $0$ or $1$. The proof of Proposition
\ref{egcond} is complete.\vspace*{-2pt}
\end{pf*}

Proposition 1.2 can be extended into a complete version of
the Erd\H{o}s--Gallai criterion for graph limits. Suppose that $W(x,y)$
is a symmetric function from $[0,1]^2$ into $[0,1]$. In
\cite{diaconisetal08}, Section 4, it is shown that the correct analog
of the degree distribution for the graph limit $W$ is the distribution
of the random variable\vspace*{-2pt}
%
%e18 ###
\begin{equation}\label{degx}
X = \int_0^1 W(U,y)\,dy,\vspace*{-2pt}
\end{equation}
where $U$ is a random variable distributed uniformly in $[0,1]$. If a
sequence of graphs converges to $W$ then the distribution of the random
variable $d_i/n$ (where $i$ is chosen uniformly from $n$ vertices and
$d_i$ is the degree of $i$) converges to $X$ in distribution. The
following result characterizes limiting degree variates.\vspace*{-2pt}

\begin{prop}
Let $X$ be a random variable with values in $[0,1]$. Let $D(x) = \sup
\{y\dvtx P(X>y) \ge x\}$. Then $X$ has the representation \eqref{degx} if
and only if for all $x\in(0,1]$\vspace*{-2pt}
\[%\label{degx2}
\int_0^x D(y)\,dy \le x^2 + \int_x^1 \min\{D(y),x\}\,dy.\vspace*{-2pt}
\]
\end{prop}

The proof is essentially as given above, approximating $W$ by a
sequence of finite graphs and using the Erd\H{o}s--Gallai criterion. We
omit further details.\vspace*{-2pt}

%s6 ###
\section{\texorpdfstring{Proof of Theorem \protect\ref{degthm} (convergence to graph limit)}%
{Proof of Theorem 1.1 (convergence to graph limit)}}\label{degthmproof}

%s6.1 ###
\subsection{Preliminary lemmas}

We need a couple of probabilistic results before we can embark on the
proof of Theorem \ref{degthm}. The first one is a simple application of
the method of bounded differences for concentration inequalities.

\begin{lmm}\label{bddiff}
Let $H$ be a finite simple graph of size $\le n$. Let $G$ be a~random
graph on $n$ vertices with independent edges. Let $t(H,G)$ be the
homomorphism density of $H$ in $G$, defined in \eqref{homdens}. Then
for any $\varepsilon>0$,
\[
\mathbb{P}\bigl(|t(H,G) - \mathbb{E}t(H,G)| > \varepsilon\bigr) \le
2e^{-C\varepsilon^2n^2},
\]
where $C$ is a constant that depends only on $H$.
\end{lmm}

\begin{pf}
The proof is a simple consequence of the bounded difference inequality
\cite{mcdiarmid89}. Note that the quantity $t(H,G)$ is a function of
the edges of $G$, considered as independent Bernoulli random variables.
When a particular edge is added or removed (i.e., the corresponding
Bernoulli variable is set equal to $1$ or $0$), $\hom(H,G)$ is altered
by at most $Cn^{|V(H)|-2}$, where $C$ is a~constant that depends only
on $H$. This is because when we fix an edge, we are fixing its two
endpoints, which leaves us the freedom of choosing the remaining
$|V(H)|-2$ vertices arbitrarily when constructing a homomorphism.

Thus, alteration of the status of an edge changes $t(H,G)$ by at most
$Cn^{-2}$. The bounded difference inequality completes the proof.
\end{pf}

The second preliminary result that we need is a kind of local limit
theorem that we need to pass from the $\bolds{\beta}$-model to graphs
with given
degree sequence.

Let $\mathbf{d}=(d_1,\ldots,d_n)$ be a valid degree sequence on a graph
of size
$n$. Let $G=(V,E)$ be a random graph on $n$ vertices labeled
$1,\ldots,n$ so that edges~$i$, $j$ are connected with probability
$p_{ij}$ satisfying $d_i = \sum_{j\neq i} p_{ij}$ and so that $\delta
\leq p_{ij} \leq1-\delta$ for some fixed $0<\delta<\frac12$. Let
$w_{ij}$ denote the indicator that $(i,j)$ is an edge in $G$. We
obtain a lower bound on the probability that $G$ has degree sequence
$\mathbf{d}$.

\begin{lmm}\label{l:lowerBound}
For any $\varepsilon>0$ and large enough $n,$ the random graph $G$ has
degree sequence $\mathbf{d}$ with probability at least
$\frac12\exp(-\log(\delta) n^{(3/2)+\varepsilon})$.
\end{lmm}

We first prove the following claim about the existence of 0--1
contingency tables. An $m\times n$ 0--1 contingency table with integer
row and column sums $r_1, \ldots, r_m$ and $c_1, \ldots, c_n$ is an
$m\times n$ matrix whose entries are $0$ or $1$ and whose $i$th row and
$j$th column sum to $r_i$ and $c_j$, respectively. Denote the
conjugate sequences as $r^*_i = \#\{r_j \dvtx r_j \geq i\}$ and $c^*_i =
\#\{c_j \dvtx c_j \geq i\}$. Let $(r_{[i]}),(c_{[i]})$ denote the order
statistics of $(r_i)$ and $(c_i)$, that is, permutations of the
sequences such that $r_{[1]}\geq r_{[2]} \geq\cdots\geq r_{[m]}$ and
$c_{[1]}\geq c_{[2]} \geq\cdots\geq c_{[m]}$.

A condition of Gale and Ryser \cite{Gale,Ryser} says that there exists
a 0--1 contingency table for row and column sums $r_1, \ldots, r_m$ and
$c_1, \ldots, c_n$ if and only if $\sum_{i=1}^m r_i = \sum_{i=1}^n c_i$
and
%
%e20 ###
%e19 ###
\begin{eqnarray}
\label{e:Gale}\sum_{i=1}^k r_{[i]} &\leq&\sum_{i=1}^k c^*_{i},\qquad 1\leq k \leq
m,
\\
\label{e:Gale2}\sum_{i=1}^k c_{[i]} &\leq&\sum_{i=1}^k r^*_{i},\qquad 1\leq k \leq
n.
\end{eqnarray}

\begin{claim}\label{c:contingencyTable}
Let $0<\delta<\frac12$ and let $(p_{ij})$ be an $m\times n$ matrix such
that $\delta\leq p_{ij} \leq1-\delta$. Suppose that $(r_i)$ and
$(c_i)$ are integer sequences satisfying the following:
\begin{itemize}
\item$\sum_{i=1}^m r_i = \sum_{i=1}^n c_i$;\vspace*{4pt}
\item$|r_i - \sum_{j=1}^n
p_{ij}| \leq\frac14\delta^2 n$ for $1\leq i \leq m$;\vspace*{4pt}
\item$|c_j -
\sum_{i=1}^m p_{ij}| \leq\frac14\delta^2 m$ for $1\leq j \leq n$.\vspace*{2pt}
\end{itemize}
Then there exists a 0--1 contingency table with row and column sums
$(r_i)$ and $(c_i)$.
\end{claim}

\begin{pf}
We establish that the Gale--Ryser conditions hold. Without loss of
generality we may assume that $r_1 \geq r_2 \geq\cdots\geq r_m$. Then
condition \eqref{e:Gale} is equivalent to
%
%e21 ###
\begin{equation}\label{e:GaleRequired}
\sum_{i=1}^k r_{i} \leq\sum_{i=1}^k c^*_{i} = \sum_{i=1}^k
\sum_{j=1}^n \mathbh{1}_{\{c_j \geq i\}} = \sum_{j=1}^n \min\{
k,c_j\}.
\end{equation}
Now
\[
\sum_{i=1}^k r_{i} \leq\sum_{i=1}^k \sum_{j=1}^n p_{ij} +
\frac14\delta^2 kn
\]
and hence,
\begin{eqnarray*}
\sum_{j=1}^n \min\{k,c_j\}
&\geq&
\sum_{j=1}^n \min\Biggl\{k,\sum_{i=1}^m p_{ij} - \frac14\delta^2 m\Biggr\}
\\
&\geq&
\sum_{i=1}^k \sum_{j=1}^n p_{ij} + \sum_{j=1}^n \min\Biggl\{k -\sum_{i=1}^k p_{ij},\sum_{i=k+1}^m p_{ij} - \frac14\delta^2 m\Biggr\}
\\
&\geq&
\sum_{i=1}^k \sum_{j=1}^n p_{ij} + \sum_{j=1}^n \min\biggl\{\delta k , (m-k)\delta- \frac14\delta^2m\biggr\}
\\
&\geq&
\sum_{i=1}^k r_{i} + n \biggl(\min\biggl\{\delta k , (m-k)\delta-\frac14\delta^2 m\biggr\} - \frac14\delta^2 k\biggr),
\end{eqnarray*}\vadjust{\goodbreak}
where we used the fact that $\delta\leq p_{ij}\leq1-\delta$. Now
$\delta k \geq\frac14\delta^2 k$ and when $1\leq k \leq
m(1-\delta+\frac14\delta^2)$,
\begin{eqnarray*}
(m-k)\delta- \tfrac14\delta^2 m
&=&
m\delta\bigl(1-\tfrac14\delta\bigr)-k\delta
\\
&\ge&
m\delta\bigl(1-\delta+ \tfrac14\delta^2\bigr)\bigl(1+\tfrac14\delta\bigr) -k\delta
\\
&\ge&
k\delta\bigl(1+\tfrac14\delta\bigr) - k\delta= \tfrac14\delta^2 k
\end{eqnarray*}
%
%(m-k)\delta- \frac14\delta^2 m \geq m (\delta- \frac14\delta^2 m) -
and hence, $n (\min\{\delta k , (m-k)\delta- \frac14\delta^2 m\} -
\frac14\delta^2 k)\geq0$. To establish equation~\eqref{e:GaleRequired}
it then suffices to consider $m(1-\delta+\frac14\delta^2)\leq k \leq
m$. In this case,
\[
c_j \leq\sum_{i=1}^m p_{ij} + \frac14\delta^2 m \leq(1-\delta) m +
\frac14\delta^2 m \leq k
\]
and so
\[
\sum_{j=1}^n \min\{k,c_j\} = \sum_{j=1}^n c_j = \sum_{i=1}^m r_{i}
\geq\sum_{i=1}^k r_{[i]}
\]
establishing \eqref{e:Gale}. Condition \eqref{e:Gale2} follows
similarly and hence, there exists a~0--1 contingency table with the
prescribed row and column sums.
\end{pf}

We are now ready to prove Lemma \ref{l:lowerBound}.

\begin{pf}
We split the $n$ vertices into subsets $A={1,2,\ldots,n-n^a}$ and
$B={n-n^a+1,\ldots,n}$ where $a =\frac12 + \varepsilon$. For $1\leq i < j
\leq|A|$, choose $w_{ij}$ according to $p_{ij}$. Let $\mathcal{G}$
denote the event that the following conditions hold:
\begin{itemize}
\item For all $i \in A$
%
%e22 ###
\begin{equation}\label{e:randomRowSum}
\biggl| \sum_{j \in A \backslash\{i\} } w_{ij} - \sum_{j \in A
\backslash\{i\} } p_{ij} \biggr| < n^{(1+\varepsilon)/2};
\end{equation}
\item That the total number of edges in the subgraph induced by $A$
satisfies
%
%e23 ###
\begin{equation}\label{e:edgeBound}
\sum_{i \in A} \biggl( d_i - \sum_{j \in A \backslash\{i\} } w_{ij} \biggr) <
\sum_{i \in B} d_{i}< \sum_{i \in A} \biggl( d_i - \sum_{j \in A
\backslash\{i\} } w_{ij}\biggr) + |B|(|B| - 1).
\end{equation}
\end{itemize}
Both conditions hold with high probability by simple applications of
Hoeffding's inequality \cite{hoeffding63}. For example, the first
follows from Hoeffding's inequality as
\[
P\biggl(\biggl| \sum_{j \in A \backslash\{i\} } (w_{ij} - \mathbb{E}w_{ij} )\biggr|
\geq
n^{(1+\varepsilon)/2}\biggr)\leq2e^{-(1/2) n^{\varepsilon}}
\]
and taking a union bound over $i\in A$.

We will show that given $\mathcal{G}$ there is always a way to add
edges between vertices in $B\times V$ so that the graph has degree
sequence $\mathbf{d}$. First we choose any assignment of the edges
$(w_{ij})_{i,j\in B}$ in $B\times B$ so that the total number of edges
equals
\[
\frac{1}{2} \biggl( \sum_{i \in B} d_{i} - \sum_{i \in A} \biggl( d_i - \sum_{j
\in
A \backslash\{i\} } w_{ij} \biggr)\biggr)
\]
which is an integer because the sum of the degrees is even and is
between 0 and $\frac12 |B|(|B|-1)$ by equation \eqref{e:edgeBound}.

It remains to assign edges between $A$ and $B$ so that the graph has
degree sequence $\mathbf{d}$. This is exactly equivalent to the
question of
finding a 0--1 contingency table with dimensions $|A| \times|B|$, row
sums $r_i = d_i - \sum_{j \in A \backslash\{i\} } w_{ij}$ for $i \in
A$ and column sums $c_i = d_i - \sum_{j \in B \backslash\{i\} }
w_{ij}$ for $i \in B$.

Condition \eqref{e:randomRowSum} guarantees that $r_i = [1+o(1)]
\sum_{j \in B} p_{ij}$ and since $|B| = o[|A|)],$ we have that $c_i =
[1+o(1)] \sum_{j \in A} p_{ij}$ uniformly in $n$. Hence, by Claim~\ref{c:contingencyTable} a 0--1 contingency table with row and column
sums $(r_i)$ and $(c_j)$ exists.

Hence, whenever the edges $(w_{ij})_{i,j\in A}$ satisfy $\mathcal{G}$
there exists at least one way to assign the other edges so that the
graph has degree sequence $\mathbf{d}$. Since any configuration
$(w_{ij})_{i\in V,j\in B}$ has probability at least $\delta^{|V||B|}$
and is independent of $\mathcal{G,}$ the probability that $G$ has the
degree sequence $\mathbf{d}$ is at least
$P(\mathcal{G})\exp(-\log(\delta)n^{1+a})$ and the result follows since
$\mathcal{G}$ holds with high probability.
\end{pf}

An alternative approach in the above lower bound could be through the
enumeration of the number of graphs of a particular degree sequence as
carried out in \cite{barvinokHartigan10}. In fact, this approach would
give a better lower bound than the one we obtain. This was brought to
our attention recently by Alexander Barvinok.
%While this is carried out in~\cite{barvinokHartigan10}, their
%enumeration is given in terms of a ``maximum entropy matrix'' which is
%not directly comparable to our bounds.

%s6.2 ###
\subsection{\texorpdfstring{Proof of Theorem \protect\ref{degthm}}{Proof of Theorem 1.1}}
Let $\mathbf{d}^n$, $G_n$ and $f$ be as in the statement of the
theorem. By
Proposition \ref{egcond} we know that $f$ has the following two
properties:
\begin{enumerate}[A.]
\item[A.] There are two constants $c_1 > 0$ and $c_2< 1$ such that $c_1\le
f(x)\le c_2$ for all $x\in[0,1]$.
\item[B.] For each $0 < b\le1$,
\[
\inf_{x\ge b}\biggl\{\int_x^1 \min\{f(y), x\}\,dy + x^2 - \int_0^x f(y)\,dy\biggr\} >0.
\]
\end{enumerate}
(The infimum is positive because the term within the brackets is a
positive continuous function of $x$.) Now fix $n$ and for each
$B\subseteq\{1,\ldots,n\}$, consider the quantity
\[
\mathcal{E}(B) := \sum_{j\notin B} \min\{d_j^n, |B|\} + |B|(|B|-1) -
\sum_{i\in B} d_i^n.
\]
Under the assumption that $d_1^n \ge d_2^n \ge\cdots\ge d_n^n$, we
claim that for each $1\le k\le n$, $\mathcal{E}(B)$ is minimized over all
subsets $B$ of size $k$ when $B = \{1,\ldots, k\}$. To prove this, take
any $B$ of size $k$. Suppose there is $a\in B$ and $b\notin B$ such
that $b < a$. Let $B' = (B \backslash\{a\})\cup\{b\}$. Then clearly,
since $d_b^n \ge d_a^n$, we have
\[
\sum_{j\notin B} \min\{d_j^n, k\} \ge\sum_{j\notin B'}
\min\{d_j^n, k\}
\]
and
\[
\sum_{i\in B} d_i^n \le\sum_{i\in B'} d_i^n.
\]
Thus, $\mathcal{E}(B) \ge\mathcal{E}(B')$, which proves the claim.
Now by the
definition of convergence of degree sequences and the fact that $f$ is
bounded and nonincreasing,
%
%e24 ###
\begin{eqnarray}\label{e:degL1Conv}
&&\Biggl| \sum_{i=1}^k \frac{1}{n}\cdot\frac{d_i^n}{n} - \int_{0}^{k/n} f(y)\,dy\Biggr|\nonumber
\\
&&\qquad\leq
\sum_{i=1}^n \biggr| \frac{1}{n}\cdot\frac{d_i^n}{n} - \int_{(i-1)/n}^{i/n} f(y)\,dy \biggr|
\\
&&\qquad\leq
\frac1n\sum_{i=1}^n \biggl|\frac{d_i^n}{n} - f \biggl(\frac{i}{n}\biggr)\biggr| +
\int_0^1
|f(x) - f(\lceil n x \rceil/n)|\,dx \to0\nonumber
\end{eqnarray}
for $1\leq k \leq n$. Similarly
\[
\sum_{j=k+1}^n \frac{1}{n}\min\biggl\{\frac{d_j^n}{n}, \frac{k}{n}\biggr\} -
\int_{k/n}^1 \min\{f(y), k/n\}\,dy \to0
\]
uniformly in $1 \leq k \leq n$ as $n \rightarrow\infty$. Hence, we
have that
for any $b\in(0,1)$,
\begin{eqnarray*}
&&
\frac{1}{n^2}\min_{B\subseteq\{1,\ldots, n\}, |B|\ge bn}\mathcal{E}(B)
\\
&&\qquad=
\min_{k\ge bn} \Biggl\{\sum_{j=k+1}^n \frac{1}{n}\min\biggl\{ \frac{d_j^n}{n},\frac{k}{n}\biggr\} + \frac{k(k-1)}{n^2} - \sum_{i=1}^k \frac{1}{n}\cdot\frac{d_i^n}{n}\Biggr\}
\\
&&\qquad\rightarrow
\inf_{x\ge b} \biggl\{\int_x^1 \min\{f(y), x\}\,dy + x^2 -\int_0^x f(y)\,dy \biggr\}\qquad\mbox{as }n \rightarrow\infty.
\end{eqnarray*}
Thus, we can apply properties A and B of the function $f$, the
definition of scaling limit of degree sequences and Lemma \ref{mainlmm}
to conclude that for all large $n$, a solution $\bolds{\beta}^n =
(\beta^n_1,
\ldots,\beta^n_n)$ to \eqref{mle} for $\mathbf{d}^n$ exists and
$|\bolds{\beta}^n|_\infty$ is uniformly bounded.

For each $n$, define a function $g_n\dvtx [0,1]\rightarrow\mathbb{R}$ as
\[
g_n(x) := \beta^n_i \qquad\mbox{if } \frac{i-1}{n} < x \leq\frac{i}{n}
\]
and let $g_n(0) := \beta^n_1$. Now fix two positive integers $m,n$ and
let
\[
N := mn.
\]
Define a vector $\mathbf{x}_0 = (x_{0,1},\ldots,x_{0,N})\in\mathbb
{R}^N$ as
follows:
\[
x_{0,i} = \beta^n_{k}\qquad \mbox{if } m(k-1) + 1 \le i\le mk.
\]
In other words,
\[
\mathbf{x}_0 = (\beta_1^n,\beta_1^n,\ldots, \beta_1^n, \beta_2^n,
\beta_2^n,\ldots,\beta_2^n, \ldots,
\beta_n^n,\beta_n^n,\ldots,\beta_n^n),
\]
where each $\beta_k^n$ is repeated $m$ times. For $\ell\geq1$ define
$\mathbf{x}_\ell=\varphi(\mathbf{x}_{\ell-1})$ as in Theorem~\ref{algo} (with
$N$ in
place of $n$). Equivalently,
%
%e25 ###
\begin{equation}\label{e:xIterate}
x_{\ell,i} - x_{\ell-1,i} = \log d_i^N - \log y_{\ell-1,i} = \log
\frac{d_i^N/N}{y_{\ell-1,i}/N},
\end{equation}
where
\[
y_{\ell,i} := \sum_{j\ne i} \frac{e^{x_{\ell,i} +
x_{\ell,j}}}{1+e^{x_{\ell,i} + x_{\ell,j}}}.
\]
Note that by definition of $y_{0,i}$ and $x_{0,i}$, if $m(k-1)+1\le
i\le mk$,
\[
y_{0,i} - m d_k^n = (m-1)\frac{e^{2\beta_k^n}}{1+e^{2\beta_k^n}}
\le
m.
\]
Consequently, if $m(k-1) + 1 \le i\le mk$,
%
%e26 ###
\begin{equation}\label{yineq}
|y_{0,i}/N - d_k^n/n| \le1/n.
\end{equation}
Hence, by equation \eqref{degconv} [similarly to \eqref{e:degL1Conv}]
it follows that
\[
\frac{1}{N}\sum_{i=1}^N |y_{0,i}/N - d_i^N/N| \leq\varepsilon_1(n)
\]
uniformly in $N$ where $\varepsilon_1(n)\to0$ as $n\to\infty$. From
\eqref{degconv}, \eqref{e:xIterate}, \eqref{yineq} (and implicitly
using the continuity of log, property A of the function $f$ and the
uniform boundedness of $|\bolds{\beta}^n|_\infty$), we see that
\[
|\mathbf{x}_0 - \mathbf{x}_1 |_1 \leq N \varepsilon_2(n)
\]
uniformly in $m$ where $\varepsilon_2(n)\to0$ as $n\to\infty$. Since
$|\bolds{\beta}^n|_\infty$ is uniformly bounded in $n$ by
Theorem~\ref{algo}, it
follows that for large enough $n,m$,
%
%e27 ###
\begin{equation}\label{e:xInftyBound}
|\mathbf{x}_\ell- \bolds{\beta}^N|_\infty\leq K \theta^\ell
\end{equation}
for some $K$ and $0<\theta<1$ independent of $n$ and $m$. Hence, for
some $K',$ also independent of $n,m$,
\[
\sup|\mathbf{x}_\ell|_\infty\leq K'.
\]
Consequently, by Lemma~\ref{l:L1Recursion} we have that
%
%e28 ###
\begin{equation}\label{e:xOneBound}
|\mathbf{x}_0 - \mathbf{x}_\ell|_1 \leq\Biggl(\sum_{i=1}^\ell
(2e^{2K'})^i \Biggr) |\mathbf{x}_0 -
\mathbf{x}_1 |_1.
\end{equation}
Combining equations~\eqref{e:xInftyBound} and~\eqref{e:xOneBound} and
using the fact that $|\mathbf{x}|_1\leq N |\mathbf{x}|_\infty$ we
have that
\[
|\mathbf{x}_0 - \bolds{\beta}^N|_1 \le|\mathbf{x}_0 - \mathbf
{x}_\ell|_1 + |\mathbf{x}_\ell-
\bolds{\beta}^N|_1
\leq\Biggl(\sum_{i=1}^\ell(2e^{2K'})^i \Biggr) N \varepsilon_2(n) + K \theta
^\ell N.
\]
Now taking $\ell=\ell(n)$ to infinity slowly enough so that
\[
\Biggl(\sum_{i=1}^\ell(2e^{2K'})^i \Biggr) \varepsilon_2(n) \to0
\]
it follows that
\[
|\mathbf{x}_0 - \bolds{\beta}^N|_1 \le N \varepsilon_3(n)
\]
uniformly in $m$ where $\varepsilon_1(n)\to0$ as $n\to\infty$. But
\[
|\mathbf{x}_0 - \bolds{\beta}^N|_1 = N\|g_n - g_N\|_1,
\]
where $\|\cdot\|_1$ is the usual $L^1$ norm on functions on $[0,1]$.
Thus,
\[
\|g_n - g_m\|_1 \le\|g_n - g_N\|_\infty+ \|g_m- g_N\|_\infty\le
\varepsilon_3(n) + \varepsilon_3(m).
\]
This shows that the sequence $\{g_n\}$ is Cauchy under the $L^1$ norm
and thus there exists a uniformly bounded function $g^*$ such that
$\|g_n - g^*\|_1\rightarrow0$. Now, for each $n$ define a function
$f_n$ as
\[
f_n(x) := \int_0^1 \frac{e^{g_n(x)+g_n(y)}}{1+e^{g_n(x)+g_n(y)}}\,dy.
\]
Now by the uniform boundedness of the $|g_n|_\infty$,
\[
\int_0^1 \biggl| f_n(x) -
\int_0^1\frac{e^{g^*(x)+g^*(y)}}{1+e^{g^*(x)+g^*(y)}}\,dy \biggr|\,dx \to0
\]
as $n\to\infty$. But from the relation between $\bolds{\beta}^n$
and $\mathbf{d}
^n$, it
is easy to see that for $x\in(0,1]$ that $f_n(x)=d^n_{\lceil n x
\rceil}/n+O(1/n)$ and hence,
\[
\lim_n \|f-f_n\|_1 \to0.
\]
It follows that
%
%e29 ###
\begin{equation}\label{e:fAlmostEverywhere}
f(x) = \int_0^1 W^*(x,y)\,dy\qquad \mbox{a.e.},
\end{equation}
where
\[
W^*(x,y) = \frac{e^{g^*(x)+g^*(y)}}{1+e^{g^*(x)+g^*(y)}}.
\]
We now adjust $g^*$ on a set of measure 0 so that
equation~\eqref{e:fAlmostEverywhere} holds for all $x$. Set
$\psi\dvtx \mathbb{R}\to(0,1)$ as
\[
\psi(z)= \int_0^1 \frac{e^{z+g^*(y)}}{1+e^{z+g^*(y)}}\,dy.
\]
By construction and since $g^*$ is uniformly bounded, it follows that
$\psi(z)$ is continuous, strictly increasing and bijective. By
equation~\eqref{e:fAlmostEverywhere} we have that
\[
f(x)=\psi(g^*(x))\qquad \mbox{a.e.}
\]
and hence, if we set
\[
g(x)=\psi^{-1}(f(x)),
\]
then $g(x)=g^*(x)$ almost everywhere. Then for all $x\in[0,1]$,
\[
f(x) = \int_0^1 W(x,y)\,dy,
\]
where
\[
W(x,y) = \frac{e^{g(x)+g(y)}}{1+e^{g(x)+g(y)}}.
\]
Moreover, by the properties of $\psi$ and $f$, we have that
$g\in{D'[0,1]}$ and its points of discontinuity are the same as $f$.

Let us now prove that $g$ is the only function in ${D'[0,1]}$ with the
above relationship with $f$. Suppose $h$ is another such function. Fix
any $n$. Define a vector $\mathbf{x}_0 = (x_{0,1},\ldots,x_{0,n})\in
\mathbb{R}^n$
as
\[
x_{0,i} := h(i/n),\qquad i=1,\ldots,n.
\]
For each $1\le i\le n$, define
\[
y_i := \sum_{j\ne i} \frac{e^{x_{0,i} + x_{0,j}}}{1+e^{x_{0,i} +
x_{0,j}}}.
\]
Then since $h\in{D'[0,1]}$,
%
%e30 ###
\begin{equation}\label{yineq2}
\sup_i |y_i/n - f(i/n)| = \sup_i \biggl|y_i/n - \int_0^1
\frac{e^{h(i/n)+h(y)}}{1+e^{h(i/n)+h(y)}} \,dy \biggr| \leq\varepsilon_4(n),
\end{equation}
where $\varepsilon_4(n) \rightarrow0$ as $n \rightarrow\infty$. Define
$\mathbf{x}_1$ in terms of
$\mathbf{x}_0$ and $\mathbf{d}^n$ as in Theorem~\ref{algo}. Then
for each $i$,
\[
x_{1,i} - x_{0,i} = \log d_i^n - \log y_i = \log\frac{d_i^n/n}{y_i/n}.
\]
From \eqref{degconv}, \eqref{yineq2} and the above identity (and
implicitly using the property A of $f$), we see that
\[
|\mathbf{x}_1 - \mathbf{x}_0|_\infty\le\varepsilon_5(n),
\]
where $\varepsilon_5(n) \rightarrow0$ as $n \rightarrow\infty$. Thus,
by Theorem \ref{algo}
we get
\[
|\mathbf{x}_0 - \bolds{\beta}^n|_\infty\le\varepsilon_6(n),
\]
where $\varepsilon_6(n) \rightarrow0$ as $n \rightarrow\infty$. This
implies that $\|h-g_n\|_1
\to0$ and hence, that $h=g$ a.e. Since we assumed both $h$ and $g$
are in ${D'[0,1]}$ this implies that~$g = h$ on $(0,1]$. To show that
$g(0)=h(0)$, observe that since $g=h$ on $(0,1]$,
\[
f(0) = \int_0^1 \frac{e^{h(0)+h(y)}}{1+e^{h(0)+h(y)}}\,dy = \int_0^1
\frac{e^{h(0)+g(y)}}{1+e^{h(0)+g(y)}}\,dy = \psi(h(0))
\]
and therefore, by the injectivity of $\psi$, $g(0)= h(0)$.

Now fix a finite simple graph $H$. Let $\bolds{\beta}^n$ be as above.
Let $G_n'$
denote a~random graph from the $\bolds{\beta}^n$-model. Let $\mathbf
{d}_n'$ be the
degree sequence of~$G'$. Then it is easy to see that conditional on the
event $\{\mathbf{d}_n' = \mathbf{d}_n\}$ the law of~$G_n'$ is the
same as that of~$G_n$.

By Lemma \ref{bddiff}, given any $\varepsilon>0$, we have that
\[
\mathbb{P}\bigl(|t(H,G_n')- \mathbb{E}t(H,G_n')| > \varepsilon\bigr) \le e^{-C_1n^2},
\]
where $C_1$ is a constant that depends only on $H$ and $\varepsilon$. By
Lemma~\ref{l:lowerBound}, we know that
\[
\mathbb{P}(\mathbf{d}_n' = \mathbf{d}_n) \ge e^{-C_2 n^{7/4}},
\]
where $C_2$ is another constant that depends only on $|\bolds{\beta
}|_\infty$.
Thus,
\begin{eqnarray*}\label{general}
\mathbb{P}\bigl(|t(H,G_n)-\mathbb{E}t(H,G_n')| > \varepsilon\bigr)
&=&
 \mathbb{P}\bigl(|t(H,G_n') - \mathbb{E}t(H,G_n')| > \varepsilon\vert\mathbf{d}_n' = \mathbf{d}_n\bigr)
 \\
&\le&
\frac{\mathbb{P}(|t(H,G_n') - \mathbb{E}t(H,G_n')| > \varepsilon)}{\mathbb{P}(\mathbf{d}_n' =\mathbf{d}_n)}
\\
&\le&
 e^{-C_3 n^2},
\end{eqnarray*}
where $C_3$ is a constant depending on $H$, $\varepsilon$ and $|\bolds
{\beta}|_\infty$.
Since $g_n \rightarrow g$, it is easy to prove that $G_n'$ converges to $W$
almost surely. From the above inequality, it follows that $G_n'$ and
$G_n$ must have the same limit almost surely. The proof of the theorem
is complete.

\section*{Acknowledgments}
The authors are indebted to Joe Blitzstein for
many helpful tips and pointers to the literature and Martin Wainwright
for the references to \cite{brown86,barndorffnielsen78,wainwrightjordan08}. We particularly thank Alexander Barvinok for
calling our attention to \cite{barvinokHartigan10} and suggesting
possible connections to our work and Svante Janson for a very careful
reading of the manuscript and pointing out numerous small errors. Last,
we thank the Associate Editor for a number of useful comments.

% imsref loaded by arune.pranskunaite, 2010-11-11 10:22:09
%

\printaddresses

\end{document}